\RequirePackage{snapshot}
\documentclass[a4paper,10pt]{article}

\usepackage[english]{babel}

\usepackage[a4paper,left=2.5cm,right=2.5cm,top=2.5cm,bottom=2.5cm]{geometry}
\usepackage{lineno}
\usepackage[final]{hyperref}
\usepackage{xcolor}
\hypersetup{
colorlinks,
citecolor=green,
linkcolor=red,
urlcolor=blue
}

\usepackage{bm}
\usepackage{authblk}
\usepackage{units}
\usepackage{graphicx}
\usepackage[toc]{appendix}
\usepackage{amsmath, amssymb}
\usepackage{cleveref} %must be after amsmath
\usepackage{authblk}
\usepackage{caption}
\usepackage{verbatim,subcaption}
\usepackage[ruled,linesnumbered]{algorithm2e}

\usepackage[numbers]{natbib}

\newcommand{\Ma}{\textrm{Ma}}
\usepackage{amsthm}

\graphicspath{ {./images/}{./images/from_giovanni/}}
\title{A reduced-order model for advection-dominated problems based on the Radon Cumulative Distribution Transform}% in $n$-dimensional space}
\author[1]{Tobias Long\thanks{Tobias.Long@nottingham.ac.uk}}
\author[1]{Robert Barnett\thanks{Robert.Barnett@nottingham.ac.uk}}
\author[2]{Richard Jefferson-Loveday\footnote{Richard.Jefferson-Loveday@nottingham.ac.uk}}
\author[3]{Giovanni Stabile\thanks{giovanni.stabile@santannapisa.it}}
\author[1]{Matteo Icardi\footnote{Matteo.Icardi@nottingham.ac.uk}}

\affil[1]{School of Mathematical Sciences, University of Nottingham, NG7 2RD, Nottingham, UK}
\affil[2]{School of Engineering, University of Nottingham, NG7 2RD, Nottingham, UK}
\affil[3]{The Biorobotics Institute, Sant'Anna School of Advanced Studies, 56025, Pontedera, Pisa, Italy }

\date{} % Leave empty to omit a date

\begin{document}

\maketitle

\begin{abstract}
    Problems with dominant advection, discontinuities, travelling features, or shape variations are widespread in computational mechanics. However, classical linear model reduction and interpolation methods typically fail to reproduce even relatively small parameter variations, making the reduced models inefficient and inaccurate. This work proposes a model order reduction approach based on the Radon-Cumulative-Distribution transform (RCDT). We demonstrate numerically that this non-linear transformation can overcome some limitations of standard proper orthogonal decomposition (POD) reconstructions and is capable of interpolating accurately some advection-dominated phenomena, although it may introduce artefacts due to the discrete forward and inverse transform. The method is tested on various test cases coming from both manufactured examples and fluid dynamics problems.
    
    % JOURNALS:
    % IJNME (International Journal for Numerical Methods in Engineering), CMAME (Computer Methods in Applied Mechanics and Engineering), CAMWA (Computers \& Mathematics with Applications),  JSC (Journal of Scientific Computing)
    
\end{abstract}

\section*{Keywords}
Reduced Order Model; Non-linear transformations, advection-dominated problems, Radon transform, Cumulative Distribution, Proper Orthogonal Decomposition

% \tableofcontents  % get rid of this once we agree on the structure

\newpage

%%%%%%%%%%%%%%%%%%%%%%%%%%%%%%
%%%%%%%%%%%%%%%%%%%%%%%%%%%%%%
\section{Introduction}
%%%%%%%%%%%%%%%%%%%%%%%%%%%%%%
%%%%%%%%%%%%%%%%%%%%%%%%%%%%%%

% Intro into importance of advective (turbulent) flows & Examples
The importance of modelling advection-dominated flows accurately and efficiently is incorporated within many scientific fields of study. Examples include the aviation ~\cite{american_institute_of_aeronautics_and_astronautics_14th_2014, kim_review_2018},
% both traditional combustion and newly emerging distributed electric propulsion (DEP) for electrical aircraft \cite{american_institute_of_aeronautics_and_astronautics_14th_2014, kim_review_2018,gohardani_synergistic_2013}. These use a configuration of electric propellers, as opposed to the typical combustion propellers in use today, allowing for more flexibility in aircraft design. In both cases
where highly turbulent and advection-dominated flow patterns are generated by aircraft propellers.
%; therefore, requiring appropriate modelling. Of particular importance to aviation, and especially DEP, has recently been noise pollution \cite{ellif_impact_2020,brentner_investigation_2018}, aiming to reduce the aeroacoustics generated by said propellers; both near-field where passengers are located, and far-field such as public on the ground. Again, this requires appropriate models of the aircraft's flow to decipher the sound frequencies produced.
Other examples of advection-dominated flow presence are the accretion disks of energy formed by black holes and certain stars \cite{b_hole_accretion_Narayan_1994};
%, b_hole_accretion_Mahadevan_1997};
aerodynamically generated noise, i.e. aeroacoustics, of helicopter rotors \cite{brentner_modeling_2003}; similarly, the near-field aerodynamics of wind turbines, close to their blades, within wind energy \cite{churchfield_advanced_2017}; and the field of turbulent flows, advective-dominance being a key contributing factor \cite{sirovich_turbulence_1987,ffowcs_williams_sound_1969,brunton_closed-loop_2015}.

% Problem: Modelling complex non-linear flows of varying geometrical configurations is demanding
The common interest in all these problems is the modelling of highly complex non-linear advection-dominated flows, aiming to predict a variety of scenarios, considering variations of the given flow that can be related to both physical or geometrical parameters. However, due to the problem's severe non-linearity and multiscale features, simulating these can be a demanding task. Since flow and transport evolutions are highly non-linear , predicting 'interim' variations -- those between two solved variational flows -- via physical interpolation of the solved flows often yields unsatisfactory results.
Moreover, the high dimensionality poses a demanding computational task, even when attempting other forms of predictive interpolation, forming significantly large data arrays to handle. As such, solving the flows in their full-dimensional form is often slow and computationally inefficient. An established methodology to combat this computational issue, producing order-reduced flows of lower dimension, is model order reduction (MOR) \cite{aroma_book,quarteroni2015reduced}.

% ROM and other approaches for advec-dom flows
MOR is widely used in computational science and engineering applications. It
is crucial in applications requiring efficient simulating systems for many
scenarios with different parameters. Example applications are system
control \cite{setia_model_2021}, uncertainty quantification
\cite{sun_non-intrusive_2021}, and optimal system design
\cite{hartmann_optimal_2014}. With model order reduction of a system, the set
of differential equations -- which describe the physical system we are
interested in -- is solved numerically in low-dimensional reduced spaces. This
contrasts with full-order models formulated in their original
high-dimensional space suitable for the problem by, for example, finite elements or finite volume methods. The
reduced space may be constructed from data in an `offline' phase and is then applied
in an `online' phase, returning an approximate solution typically in a far
shorter time than solving for the full, non-reduced model. This can dramatically
increase the speed at which many different design parameters can be tested. However, advection-dominated phenomena have seen limited use of MOR so far. This is
because breaking the Kolmogorov $n$-width barrier and 
going beyond linear MOR strategies remains a challenging problem to solve \cite{ohlberger_rave_2016, greif_urban_2019}. 
Advection-dominated problems,  in fact, exhibit a slow decay in Kolmogorov $n$-width, thus
rendering the reduced order model (ROM) construction inefficient
\cite{peherstorfer_2020}. Current approaches to overcome this problem include
the use of nonlinear compression strategies (eg. autoencoders, convolutional
autoencoders) \cite{Lee2020,RomorStabileRozza2022}, model reduction on metric spaces \cite{Ehrlacher_2020}, piecewise linear subspaces
(involving clustering into the parameter space or solution space)
\cite{Amsallem2012a}, and nonlinear transformation of the solution manifold in
order to subsequently apply standard compression strategies. This last approach can
 be divided into a few distinct approaches: where the transformation from
the physical problem is known \cite{Reiss2018}, where the transformation is
computed by solving an optimal transport problem \cite{Bernard2018}, and finally
where the transformation is learnt with a neural network
\cite{PapapiccoDemoGirfoglioStabileRozza2021}.

% Our approach and its benefits
The present work focuses on a specific non-linear transformation, called Radon CDT (RCDT), which combines the Radon transform and the recently proposed cumulative distribution transform (CDT) \cite{kolouri2015radon,park_cumulative_2018}. This approach converts two- and three-dimensional fields, e.g., a travelling wave, to a one-dimensional series (parametrised by the Radon angles). These are then interpreted as one-dimensional probability density functions and processed through the CDT. Through the conjunction of these two transforms, two and three-dimensional transport maps become linear in RCDT space  \cite{ren_model_2021}.
In image processing, this is referred to as the linear separability property of the RCDT, i.e.,  the property of snapshots to be linearly dependent, and, therefore, compressed and interpolated exactly with linear MOR methods, such as proper orthogonal decomposition (POD).
 Since transport becomes separable in the RCDT space, interpolation can be done through standard approaches, such as linear interpolation. The result is a reduced-order model (ROM) in a non-linear space whose outputs can be transformed back into physical space. This combination has been introduced in \cite{ren_model_2021} for the first time (to the best of our knowledge). Starting from this work, in our article, we introduced the following novel aspects:
\begin{itemize}
    \item In \cite{ren_model_2021} the generated reduced order space has not been used to create a parametric ROM, and the discussion on the extrapolation properties of such a ROM is missing. In fact, in \cite{ren_model_2021}, the RADON-CDT space is only used to compress and reconstruct training snapshots and not to create a ``real'' reduced order model.
    \item In \cite{ren_model_2021}, the methodology is tested on a simple one-dimensional example or manufactured test cases. In our article, we test the methods on several test cases coming from both manufactured examples and physics-based simulations. 
    \item In this paper, we provide a comprehensive discussion of the advantages and drawbacks of the methodology. This is missing in \cite{ren_model_2021}.
\end{itemize}

% What this paper entails
In this work, we utilise the peculiar properties of the RCDT to capture geometric and spatial variations within a parameterised input and use this to produce an approximate solution for system parameters in a model order reduction methodology. Initially, we investigate the properties of the RCDT with simplified test cases to gauge the strengths and weaknesses of the potential use of the transform in the ROM and CFD communities. The RCDT is then applied to a POD-based ROM workflow and later tested on a number of computational fluid dynamics (CFD) data sets. This allows for the preservation of the flow features when transformed between spaces, alongside the accuracy of ROM's flow reconstruction at reduced order. We finally introduce interpolation and study the error in predicting flows, giving an initial qualitative gauge of RCDTs' applicability to fluid dynamics and advection-dominated problems.

% Used/Implemented: RCDT workflow + ROM/POD => method of snapshots
Both the implementation of the RCDT and ROM workflows have been written in \textit{Python 3.9.7}, making use of two packages, \textit{PyTransKit} \cite{abu_hasnat_mohammad_rubaiyat_pytranskit_2022} and \textit{EZyRB} \cite{demo_ezyrb_2018}; implementing the discretised form of RCDT -- with subsequent forward/inverse transforms -- and model order reduction functionality, respectively. For the ROM side, i.e. \textit{EZyRB}, model reduction is approached using proper orthogonal decomposition (POD), see \cite{demo2019marine,tezzele2019marine}, for example, in applying \textit{EZyRB} towards shape optimization problems. All the code developed for the preparation of this article is available open-source \cite{rcdt-rom-lib}. Specifically, the singular value decomposition (SVD) -- discussed more in \cref{POD} -- is used to determine the POD modes for the reduced-order model. SVD is not the only way to compute the POD, though an alternative approach is given by the method of snapshots \cite{sirovich_turbulence_1987,wang2016approximate}. Three distinct workflows have been implemented: the RCDT transform upon a single snapshot image/flow followed by the inverse transformation to observe the \textit{intrinsic} error induced by the discretisation and implementation of the non-linear transform; the RCDT-POD reconstruction/projection error to evaluate the effect of the non-linear transformation in the POD modes; and, the complete RCDT-POD ROM workflow consisting of the RCDT-POD on a series of snapshots, and the subsequent interpolation (with respect to time or other parameters) to predict unseen scenarios.

% How we analyse i.e. errors
To analyse the RCDT, the POD, and the RCDT-POD ROM's applicability towards modelling advection-dominated flows, we consider the relative errors introduced by these transformations by performing several test cases to observe each of their magnitudes. We distinguish between different sources of errors, namely the \textit{intrinsic} error given by the RCDT discretisation and its inverse; the \textit{reconstruction/projection} error of RCDT-POD, influenced by the number of POD modes when reconstructing the flow after order-reduction; and RCDT-POD ROM \textit{interpolation} error, introduced from predicting between "known" flow configurations to 
gauge the capability of predicting, via interpolation in their respective space, along a given parameter, e.g., time, distance, 'thickness', etc.

% Structure
The paper is organised as follows: RCDT, POD and interpolation methodology are introduced in \cref{Method}, and their numerical implementation with simple test cases is reported in \cref{RCDT Numerical}. To observe and test the errors described above and assess the capabilities of RCDT in the context of MOR, we consider a number of images and flow case studies in \cref{sec::RCDTverify} and \cref{Applications}. Test cases to quantify specific errors can be found in \cref{sec::RCDTverify}. These include discontinuous images of a unit circle and a circular ring and continuous Gaussian functions to observe the {intrinsic} error in the RCDT workflow, and the flow field given by a twin jet at different separation widths to test the interpolation in RCDT space compared to the physical.
In \cref{Applications}, instead, we focus on the complete MOR procedure starting with a simple moving Gaussian distribution, transformed into RCDT space and order-reduced using POD, compared alongside 'standard' POD in physical space. We then test our workflow for a multi-phase fluid wave and the flow around an airfoil using high-resolution CFD data. Final discussions and future work directions are then reported in \cref{Conclusion}.

%%%%%%%%%%%%%%%%%%%%%%%%%%%%%%
%%%%%%%%%%%%%%%%%%%%%%%%%%%%%%
%%%%%%%%%%%%%%%%%%%%%%%%%%%%%%
%%%%%%%%%%%%%%%%%%%%%%%%%%%%%%
\section{Methodology} \label{Method}
%%%%%%%%%%%%%%%%%%%%%%%%%%%%%%
%%%%%%%%%%%%%%%%%%%%%%%%%%%%%%
%%%%%%%%%%%%%%%%%%%%%%%%%%%%%%
%%%%%%%%%%%%%%%%%%%%%%%%%%%%%%

In this section, we provide definitions of the Radon and Cumulative Distribution transformations and combine the two to reduce and process high-dimensional functions as a series of one-dimensional probability density functions. We also introduce the proper orthogonal decomposition (POD) method as a compression and interpolation tool.

%%%%%%%%%%%%%%%%%%%%%%%%%%%%%%
%%%%%%%%%%%%%%%%%%%%%%%%%%%%%%
\subsection{Radon Transform}\label{sec:RadonDef}
%%%%%%%%%%%%%%%%%%%%%%%%%%%%%%
%%%%%%%%%%%%%%%%%%%%%%%%%%%%%%

The Radon transform is a widely used image transform in fields such as medical imaging and optics. Its inverse transform (iRadon) can be used to reconstruct images from probe-measured distributions, see \cite{deans_radon_1983, natterer_mathematics_2001}. The Radon transform $\mathcal{R}f:\mathcal{S}^{n-1}\times \mathbb{R} \to \mathbb{R}$ of a function $f:\mathbb{R}^n \rightarrow \mathbb{R}$ can be defined as 

\begin{equation}
    \mathcal{R}f(\bm{\xi},p) = \int_{\mathbf{x}\cdot \bm{\xi} = p} f(\mathbf{x})\mathrm{d}m(\mathbf{x}),
\end{equation}
\\
where $\mathrm{d}m$ is the Euclidean measure over the hyperplane. For any fixed pair $(\bm{\xi},p) \in \mathcal{S}^{n-1}\times \mathbb{R}$ the set $\{ \bm{x}\in\mathbb{R}^n | p = \bm{x}\cdot \bm{\xi} \}$ defines a hyperplane. Therefore, the transform is an integration of the function over this hyperplane. In $\mathbb{R}^2$, the Radon transform integrates a function $f(x,y)$ over a series of lines and can be written explicitly as

\begin{equation}
    \mathcal{R}f(\bm{\xi},p) = \int_{-\infty}^{\infty}\int_{-\infty}^{\,\infty}f(x,y)\,\delta(p-x\,\mathrm{cos}\,\theta - y\,\mathrm{sin}\,\theta)\,\mathrm{d}x\,\mathrm{d}y.
\end{equation}
\\
Here $p$ denotes the oriented distance of the line $l_{\bm{\xi},p}$ to the origin, and $\theta$ is the angle as defined in polar coordinates that identifies points on the circle through the following definition $\bm{\xi}:=(\cos \theta,\sin \theta)$. Here $\theta$ lies in the range $0 \leq \theta < \pi $ due to the symmetry of line integrals for a 2D function in polar coordinates.

%%%%%%%%%%%%%%%%%%%%%%%%%%%%%%
%%%%%%%%%%%%%%%%%%%%%%%%%%%%%%
\subsection{Cumulative Distribution Transform}\label{sec:CDTDef}
%%%%%%%%%%%%%%%%%%%%%%%%%%%%%%
%%%%%%%%%%%%%%%%%%%%%%%%%%%%%%

The Cumulative Distribution Transform (CDT)\cite{park_cumulative_2018} is a non-linear transform for one-dimensional probability density functions (PDF). The CDT has been developed to aid in object recognition problems, rendering certain types of classification problems linearly separable in the transform space (data classification examples include hand gestures, accelerometer data \cite{park_cumulative_2018}, wave signals \cite{aldroubi_signed_2022} and optimal mass transport signal processing, see \cite{optMassTransportOverview&RCDT_Kolouri2017}). The CDT has been developed from work on linear optimal transport methods; in contrast to linear transform methods such as Fourier and Wavelet, which consider intensity variation at fixed coordinates, the CDT considers the location of intensity variations within a signal.

Let us consider a PDF $f:Y\rightarrow \mathbb{R}$ and a reference PDF $r:X\rightarrow \mathbb{R}$ where $X=[a,b]\in\mathbb{R}$ and $Y=[c,d]\in\mathbb{R}$. 
In the following, we will consider the case of the reference density to be uniform $r = U(0,1)$, i.e., $a=0$, $b=1$.

The CDT is defined as:
\begin{equation}
    \hat{f}(x) = (\lambda(x)-x)\sqrt{r(x)},     
\end{equation}
where $\lambda(x):X\to Y$  satisfies the relation:
\begin{equation}
    \int_{\inf{(Y)}}^{{\lambda}(x)}f(\tau)\,\mathrm{d}\tau = \int_{\inf{(X)}}^{x}r(\tau)\,\mathrm{d}\tau.
\end{equation}

If we introduce the cumulative distribution functions (CDF) for $f$ and $r$, $J_f:Y\to[0,1]$ and $J_r:X\to[0,1]$ respectively, defined as:
\begin{equation}
    J_f(y) = \int_{\inf{Y}}^y f(\tau) \text{d}\tau, \quad J_r(x) = \int_{\inf{X}}^x r(\tau) \text{d}\tau\,,
\end{equation}
we have the identity:
\begin{equation}
    J_r(x) = J_f(\lambda(x)),
\end{equation}
and
\begin{equation}
    r(x) = \lambda'(x)f(\lambda(x)).
\end{equation}

If $f$ is strictly positive and continuous, due to the properties of the CDF and the positivity of $f$, we have that $J_f$ is a strictly increasing (invertible) function, differentiable and with a differentiable inverse.
Under these conditions, $\lambda$ is differentiable and we can write explicitly:
\begin{equation}
    \lambda(x) = J_f^{-1}(J_r(x))\,,
\end{equation}
\begin{equation}
    \lambda'(x) = \frac{d}{dx}J_f^{-1}(J_r(x))
    = \frac{r(x)}{f(J_f^{-1}(J_r(x)))}.
\end{equation}

The inverse Cumulative distribution Transform of $\hat{f}$ is then defined as:
\begin{equation}
    f(y) = \frac{d}{dx}J_r(\lambda^{-1}(x)) = (\lambda^{-1})'(x)r(\lambda^{-1}(x))= \frac{r(\lambda^{-1}(y))}{\lambda'(\lambda^{-1}(y))}.
\end{equation}

For the case of distributions that are non-strictly positive and discontinuous in a finite number of points, the CDF is piecewise invertible and differentiable and the above definitions are valid in a piecewise manner.

The CDT carries important properties such as that translation and scaling in the input variable $x$ result, respectively in a constant  and linear lifting of the output. Namely, if $g(y)= \frac{f(y/\sigma - \mu)}{\sigma}$ is the scaled and shifted distribution then, combining the results in \cite{park_cumulative_2018}, it can be shown that:
\begin{equation}
    \hat{g}(x) = \hat{f}(x)
    + \left(\mu + x(\sigma-1)\right) \sqrt{r(x)} . 
\end{equation}
This property is the key to a good approximation of advection-dominated and transport problems. To better clarify this point we report here a prototypical example to show some properties of the CDT. Let us consider a solution manifold given by the translation of a density function $I_1(x)$ by $\mu$:

$$I_\mu (x) = I_1(x-\mu) = \frac{1}{\sqrt{2\pi}}e^{-(x-\mu)^2/2}.$$

Such a case would have a really slow decay of the Kolmogorov $n$-width, since it is a pure translation of a reference distribution. This example could represent a 1D traveling wave moving without deformations in a 1D domain. The corresponding CDT, $\hat{I}_\mu:[0,1] \to \mathbb{R}$, computed for $I_\mu$ with respect to the uniform reference density $I_0:[0,1]\to \mathbb{R}$ is given by:
$$\hat{I}_\mu (x) = \hat{I}_1(x) + \mu,$$

which is the translation constant $\mu$ plus the CDT of the $I_1$ distribution. It is clear that now the Kolmogorov $n$-width decay in the transformed space is much faster and would actually require only one POD mode to approximate the entire solution manifold. For more details on the derivation and on the mathematical proofs we refer to \cite{park_cumulative_2018}.

We can further extend the CDT for arbitrary positive functions by applying a normalisation before applying the CDT and de-normalise the result after the inversion, and for non-positive functions by performing the CDT transform separately on the positive and negative parts.
Defined as such, the CDT becomes an invertible nonlinear transform, transforming continuous probability densities to piecewise differentiable functions.  A consequence, however, of this highly nonlinear transform is that its discretisation (we refer to \cite{park_cumulative_2018} for the description of the discretised CDT) introduces an \textit{intrinsic} error when inverting back to physical space.

\subsubsection{One dimensional numerical example using the Cumulative Distribution Transform}
Here we show the application of the Cumulative Distribution Transform to a manufactured example that resembles the main features of advection-dominated problems. In particular, given the parameter space $\mathbb{P} = \left\{ \mu \in [0.25,0.75] \right\}$, we consider the following parametrised Gaussian function $f:\mathbb{R}\times \mathbb{P} \to \mathbb{R}^+$:

$$f(x;\mu) = \frac{1}{0.1 \cdot \sqrt{2\pi}} e^{-\frac{1}{2}\left(\frac{x-\mu}{0.1}\right)^2}.$$

\begin{figure}[htbp]
    \centering
\includegraphics[width=0.48\textwidth]{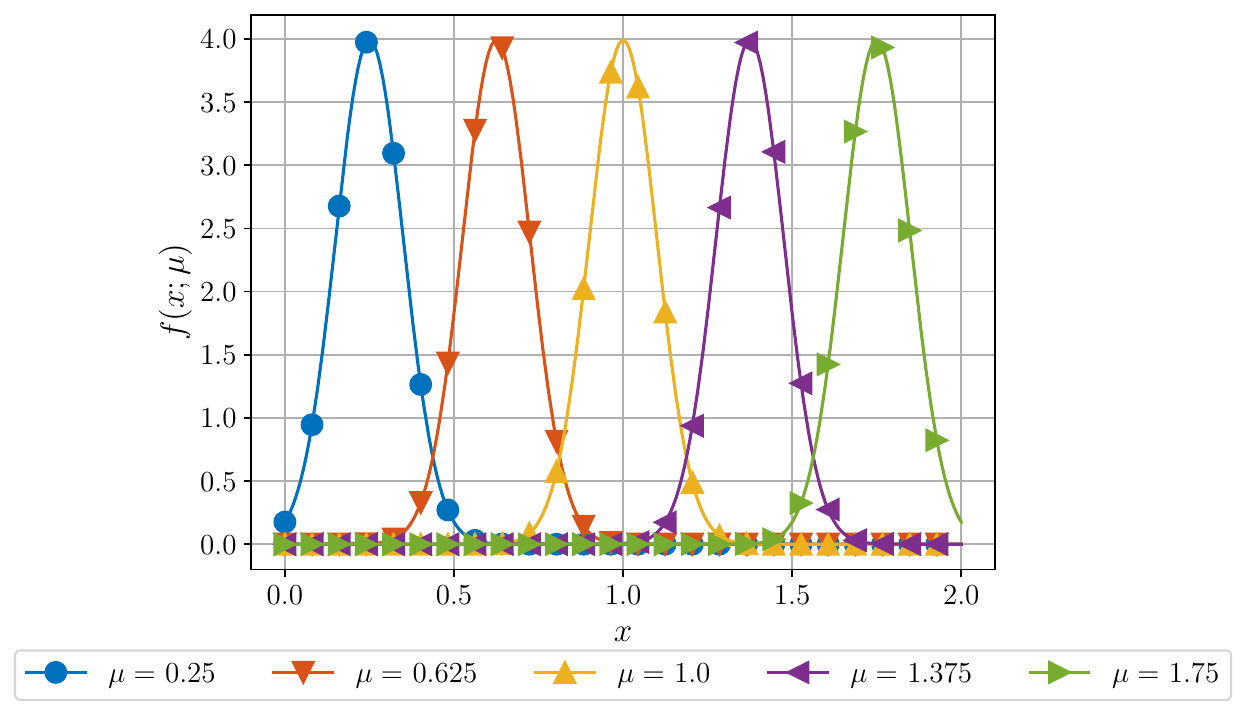}
\includegraphics[width=0.48\textwidth]{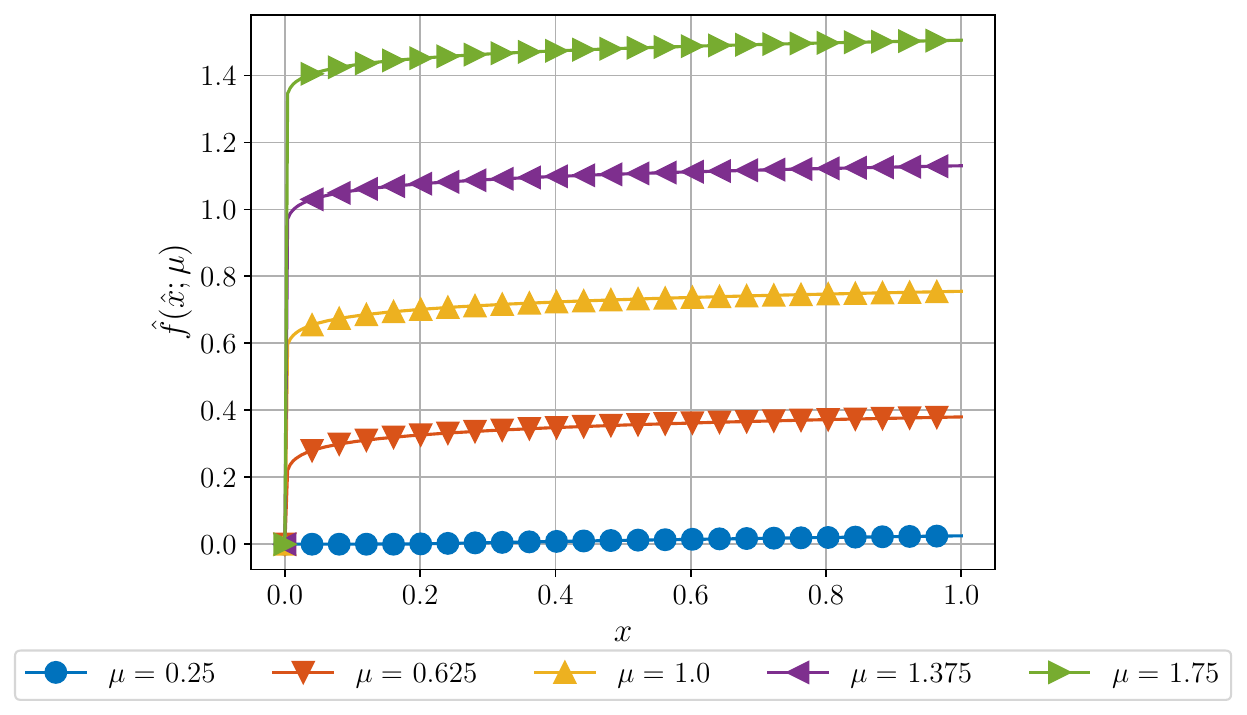}    
    \caption{Snapshots of the manufactured example for $5$ different values on the input parameter in the physical space (on the left) and in the CDT space (on the right).}\label{fig:1d_snaps}
\end{figure}

\begin{figure}[htbp]
    \centering
    \includegraphics[width=0.4\textwidth]{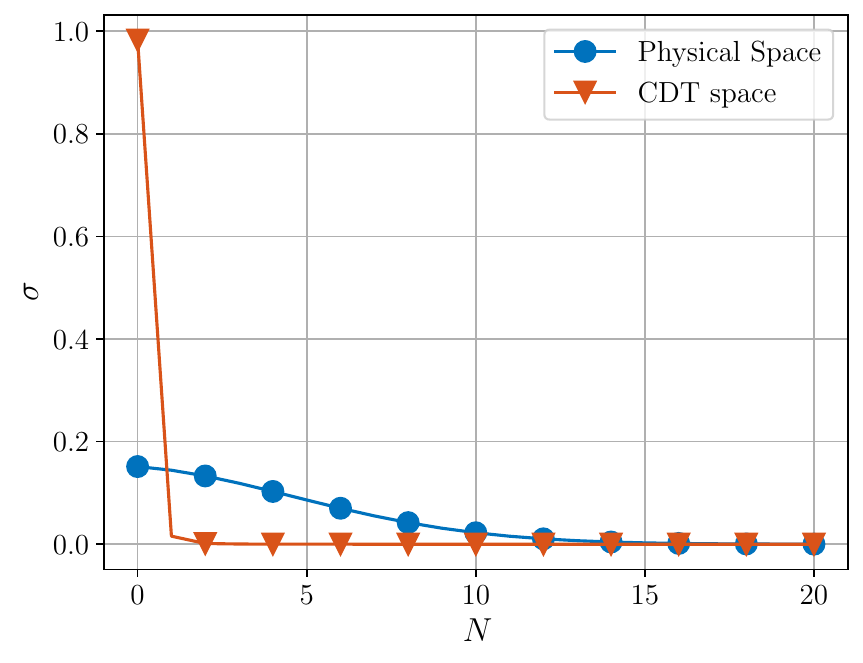}
    \includegraphics[width=0.4\textwidth]{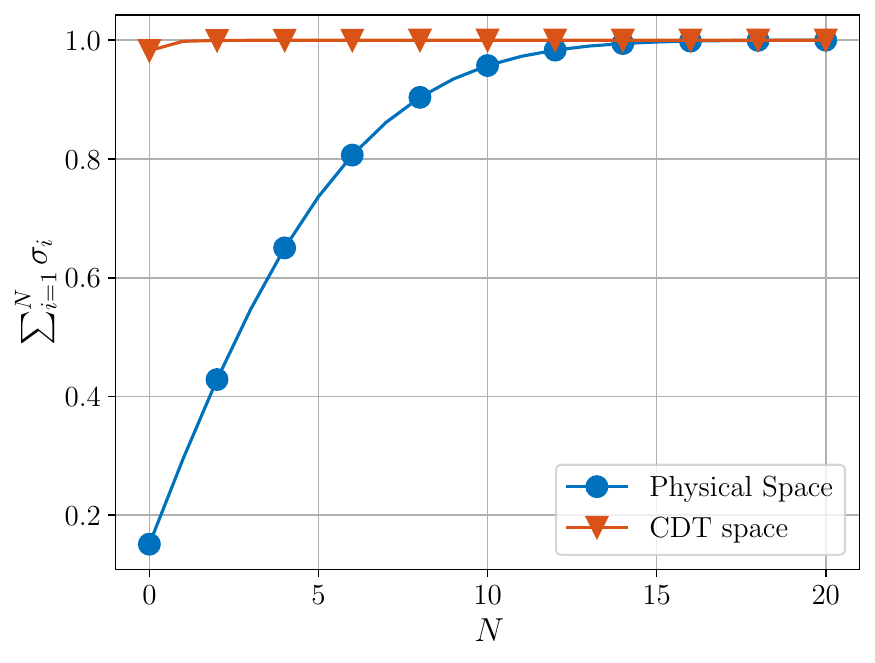}
    \caption{Eigenvalue decay (left) and cumulative eigenvalues (right) of the POD computed for the manufactured example. Results are plotted both in the physical and in the CDT space.}\label{fig:1d_pod}
\end{figure}

This function wants to mimic a travelling wave that is moving along the $x$-direction. In \cref{fig:1d_snaps} we show the function $f$ for five different values of the input parameter both in the physical space (left) and in the CDT space (on the left). In \cref{fig:1d_pod} we show the comparison of the eigenvalue decay (on the left) and of the cumulative eigenvalues (right) both in the physical and in the CDT space. In the example, the POD has been computed taking $20$ snapshots calculated considering $20$ equispaced parameter values taken inside the parameter space $\mathbb{P}$.


%%%%%%%%%%%%%%%%%%%%%%%%%%%%%%
%%%%%%%%%%%%%%%%%%%%%%%%%%%%%%
\subsection{Radon Cumulative Distribution Transform}\label{sec:RadonCDTDef}
%%%%%%%%%%%%%%%%%%%%%%%%%%%%%%
%%%%%%%%%%%%%%%%%%%%%%%%%%%%%%

The Radon Cumulative Distribution Transform (RCDT) \cite{park_cumulative_2018,shifat-e-rabbi_radon_2021,ren_model_2021} is a combination of the cumulative distribution transform (CDT) and the Radon transform. The Radon transform is used to reduce high-dimensional functions to a series of 1-D functions; these are then normalised, split into positive and negative parts, and processed with the CDT. This allows us to expand the use of one-dimensional CDT to higher dimensions, such as two or three-dimensional flow velocity data. In this work, we utilise the RCDT to transform 2D data in uniformly-spaced grids, which can be interpreted as images. For 3D data sets, 2D slices along one dimension of the data can be used with the 2D RCDT procedure or, alternatively, the 3D Radon transform algorithm can be applied.

The RCDT applied to a 2D function $f:[a,b]\times[c,d]\to\mathbb{R}$ is therefore defined as $\tilde{f}:[0,\pi]\times[s_1,s_2] \to \mathbb{R}$  by applying the CDT for each projection angle in the Radon transform, which results in
\begin{equation}
    \tilde{f}(\theta,s) = \left(\lambda(\theta,s) - s\right)\sqrt{r(s)},
\end{equation}
where
\begin{equation}
    \int_{s_1}^{\lambda(s,\theta)} \mathcal{R}f((\cos\theta,\sin\theta),s) \, \mathrm{d}s' = \int_{s_1}^s r(s') \, \mathrm{d}s'\,,
\end{equation}
and $r(s)=U(s_1,s_2)$ is the reference uniform density on the interval $[s_1,s_2]$ representing the projection of the support of the Radon transform in the second variable.
Due to the composition of two invertible transforms, the RCDT transform is invertible.
For the interested reader, details of the computational algorithms for the discrete versions of the CDT and combined RCDT that are used here are laid out in \cite{park_cumulative_2018} and \cite{park_cumulative_2018} respectively. In previous applications within image processing, the RCDT has been shown to capture scaling and transport behaviours well, allowing interpolation between geometrical features within data to great success. Note that this feature arises within the CDT rather than the Radon transform. 

The Radon transformation is used as a preprocessing tool to transform a 2D or 3D image (snapshot) into one-dimensional functions which can be processed with the CDT transformation. We refer to \cite{park_cumulative_2018} for more details on the properties of the RCDT applied to image processing. 

For example, the RCDT was successfully used for interpolation between pictures of human faces \cite{park_cumulative_2018} and used to increase the efficiency of image classification with neural networks \cite{shifat-e-rabbi_radon_2021}.
%In this work, we make use of the \textit{PyTransKit} package to carry out the RCDT, which has been developed by Rubaiyat \textit{et al.} \cite{abu_hasnat_mohammad_rubaiyat_pytranskit_2022}.
A caveat, however, that comes from the RCDT is the introduction of what we define as an \textit{intrinsic} error due to the forward and inverse discrete transforms between the original and RCDT space. Later in \cref{RCDT Numerical} this error is visible particularly near the image boundaries when performing RCDT, and its inverse compared to the original image.

\subsection{A remark on the \textit{intrinsic} error}
The RCDT combines two invertible transformations. Although this procedure, at the continuous level, is exact, the numerical implementation requires a series of approximations that introduce different sources of errors:
\begin{itemize}
    \item The number of considered angles in the Radon transform is finite. Therefore, there is an approximation in the computation of the Radon transform that is associated with a finite number of integration lines. 
    \item Functions that we consider are not continuous but expressed on a finite number of points (pixels). Therefore, a proper interpolation is required. Such interpolation introduces a certain amount of error. 
    \item In our approach, the computation of the inverse Radon transform is based on the Filtered BackProjection (FBP) algorithm \cite{Dudgeon1995-vh}. The mathematical foundation of the filtered back projection is the Fourier slice theorem. It uses the Fourier transform of the projection and interpolation in Fourier space to obtain the 2D Fourier transform of the image, which is then inverted to form the reconstructed image. The filtered back projection is among the fastest methods of performing the inverse Radon transform. The only tunable parameter for the FBP is the filter, which is applied to the Fourier-transformed projections. It may be used to suppress high-frequency noise in the reconstruction. This algorithm introduces an error, and the reconstructed image is different with respect to the original one. 
    \item The CDT is defined for continuous-time functions 
    in a contiguous, finite domain. To compute the CDT for a discrete function, we need a way of estimating its cumulative function at any
    arbitrary coordinate. We do so via interpolation using B-splines \cite{Unser1993} of degree zero, guaranteeing that the reconstructed functions are always positive. This interpolation procedure introduces an additional source of error. 
\end{itemize}

In the following, we denote as \textit{intrinsic} all the errors introduced by the discretisation and the numerical implementation of the RCDT and its inverse, to distinguish them from the errors introduced by the model reduction compression and interpolation procedures.
It is beyond the scope of this article to delve into details of the computational algorithms used to compute the Radon and the CDT transforms and the associated errors introduced by the discretisation. We refer interested readers to \cite{Beckmann2020,park_cumulative_2018} for more details. 

%%%%%%%%%%%%%%%%%%%%%%%%%%%%%%
%%%%%%%%%%%%%%%%%%%%%%%%%%%%%%
\subsection{Model Order Reduction with Proper Orthogonal Decomposition} \label{POD}
%%%%%%%%%%%%%%%%%%%%%%%%%%%%%%
%%%%%%%%%%%%%%%%%%%%%%%%%%%%%%

Proper Orthogonal Decomposition (POD) is the most widely used technique to
compress the solution manifold of a variety of problems like the unsteady
Navier-Stokes equations and it has been applied ubiquitously in MOR research
during the past decades. The method dates back to the work of \textit{Pearson}
\cite{pearson_liii_1901}, with POD first being applied to turbulent flows by
\textit{Lumley} in 1967 \cite{lumley_structure_1967}. POD is closely related to
methods in other areas of mathematics, e.g., principal component analysis (PCA)
in statistical analysis \cite{jolliffe_principal_2016} and the Karhunen–Loève
expansion in stochastic modelling \cite{cho_karhunenloeve_2013}. POD is a
technique for computing an orthonormal reduced basis for a given set of
experimental, theoretical or computational data. Specifically, the POD basis is obtained
by performing the singular value decomposition of the snapshot matrix, 
assembled with sampled data at distinct time or parameter instances (`snapshots') since POD is not only for time-dependent but also for parameter-dependent problems. The
computation is described for the example solution $\mathbf{u}(\mathbf{x},t)$,
for $\mathbf{x}\in\mathbb{R}^3$ and $0\leq t\leq T$. The numerical implementation for the computation of the POD in this work is handled by the python package \textit{EZyRB} \cite{demo_ezyrb_2018}.

SVD is a way to obtain the POD of a dynamical system when given a sequence of data snapshots at various time instances. First, we seek an approximation to the solution $\mathbf{u}(\mathbf{x},t)$ as a linear combination of spatial modes $\Psi_i(\mathbf{x})$ and temporal coefficients $a_i(t)$ for $i=1,\dots,N_r$. That is,
\begin{align}
    \mathbf{u}(\mathbf{x},t) \approx \mathbf{u}_r(\mathbf{x},t) = \sum_{i=1}^{N_r} a_i(t) \Psi_i(\mathbf{x}).
\end{align}
Here, $N_r \geq 1$ is the number of modes. This is reasonable for a fluid flow that can be approximated as a stochastic and stationary process in time and ergodic. The spatial modes are orthogonal, that is to say $\langle \Psi_i, \Psi_j \rangle = 0$ for $i \neq j$ where $\langle \cdot, \cdot\rangle$ denotes the $L^2$ inner product\footnote{Here we denote with $L^2$ inner product the standard Euclidean inner product instead of the $L^2$ (function space) inner product. }. The time coefficients can be calculated using various methods, the most common being via Galerkin or least square Petrov-Galerkin projection of the original system onto the spatial modes, where the resulting ROM is then termed a \textit{(Petrov-)Galerkin ROM} or \textit{POD-(Petrov-)Galerkin ROM} \cite{Holmes2012,StabileBallarinZuccarinoRozza2019,BustoStabileRozzaCendon2019}. These approaches fall into the category of intrusive ROMs. Other methods for developing a ROM from POD are \textit{POD with interpolation (POD-I)} \cite{BuiThanh2003}, \textit{POD with Gaussian process regression (POD-GPR)} \cite{Guo2018} and \textit{POD with neural networks (POD-NN)} \cite{Hesthaven2018}. These approaches fall into the category of non-intrusive ROMs \cite{chap9_aroma_book}. Intrusive methods exploit the discretised equations and project them onto the space spanned by the POD modes to obtain a lower-dimensional system of ODEs. This technique, exploiting the governing equations, is naturally \emph{physics-based} and, usually, has better extrapolation properties. However, it requires a larger implementation effort and access to the discretised differential operators assembled by the full-order model (hence the name intrusive). Moreover, in the case of non-linear/non-affine problems, the speed-up that can be achieved is usually smaller compared to non-intrusive approaches \cite{StabileZancanaroRozza2020}. On the other hand, non-intrusive methods are purely data-driven and can guarantee a speed-up also in non-linear/non-affine cases. In this article, we will focus our attention only on non-intrusive methods. The modes $\Psi_i(\mathbf{x})$ of the decomposition are calculated from snapshots of the data. The data snapshots are taken as semi-discretised state solutions of the system, computed at different points in time or at different parameter values, denoted by $\mathbf{u}_1, \mathbf{u}_2, ... , \mathbf{u}_n$ where $\mathbf{u}_j = \mathbf{u}(t_j,\mathbf{p}_j) \in \mathbb{R}^N$ denotes the $j^{th}$ snapshot at time $t_j$ and parameter values $\mathbf{p}_j$, with $n$ being the total number of snapshots. Define a snapshot matrix $\mathcal{U} \in \mathbb{R}^{N \times n}$ constructed with columns as each snapshot $\mathbf{x}_j$. The singular value decomposition of the snapshot matrix $\mathbf{\mathcal{U}}$ is then: %
\begin{align}\label{eqn:SVD}
    \mathbf{\mathcal{U}} = \mathbf{U}\mathbf{\Sigma}\mathbf{V}^T \approx \sum_{j=1}^{N_r} \sigma_j \mathbf{u}_j \mathbf{v}_j^T\,,
\end{align}
where the columns of $\mathbf{U} \in \mathbb{R}^{N \times N}$ and $\mathbf{V} \in \mathbb{R}^{n \times n}$ are the left ($\mathbf{u}_j$) and the right ($\mathbf{v}_j$) singular vectors of matrix $\mathbf{\mathcal{U}}$ respectively. The singular values of $\mathbf{\mathcal{U}}$, $\sigma_1 \geq \sigma_2 \geq ... \geq \sigma_M \geq 0$, where $M=\min(n,N)$, are found along the diagonal of the rectangular matrix $\Sigma \in \mathbb{R}^{N \times n}$. The POD basis vectors can then be chosen as the left singular vectors corresponding to the $N_r$ largest singular values. The $N_r$ singular values are generally chosen based on the ratio of the `energy' contained in the $N_r$ modes and the total `energy' in all $n$ modes, that is, the ratio given by:
\begin{align}
    \frac{\sum_{i=1}^{N_r} \sigma_i^2}{\sum_{i=1}^n \sigma_i^2} = \kappa.
\end{align}
The POD can be applied to both the time domain in time-dependent problems or the parameter domain in parametric flow problems; in both cases, the method remains the same, and the snapshots are sampled in either time or parameter space.

An obvious consequence of taking the reduced $N_r$ modes is the subsequent lost 'energy' of the $n-N_r$ modes, resulting in a \textit{reconstruction} error of the flow or moving object geometry. Taking a sufficiently large number of modes reconciles these errors but, as seen in the test cases of \cref{Applications}, this highly depends on the pre-processing and transformations done on the original snapshots.

%%%%%%%%%%%%%%%%%%%%%%%%%%%%%%
%%%%%%%%%%%%%%%%%%%%%%%%%%%%%%
%%%%%%%%%%%%%%%%%%%%%%%%%%%%%%
%%%%%%%%%%%%%%%%%%%%%%%%%%%%%%
\section{RCDT implementation and applications} \label{RCDT Numerical}
%%%%%%%%%%%%%%%%%%%%%%%%%%%%%%
%%%%%%%%%%%%%%%%%%%%%%%%%%%%%%
%%%%%%%%%%%%%%%%%%%%%%%%%%%%%%
%%%%%%%%%%%%%%%%%%%%%%%%%%%%%%

This section focuses on the algorithmic procedure of RCDT and testing its consistency in image capturing and interpolation. The implemented steps of the RCDT process have been written using \textit{Python 3.9.7}. Construction of the discrete RCDT image is done via use of the package \textit{PyTransKit} \cite{abu_hasnat_mohammad_rubaiyat_pytranskit_2022}, used for example in \cite{shifat-e-rabbi_radon_2021}, containing the RCDT class for forward and inverse transforms of two-dimensional images, defined as a data array, assuming the image is normalised.

As discussed in \cref{Method}, both CDT and RCDT have been developed for distributions, respectively, in 1D and 2D. One of the aims of this work is to extend these procedures to generic images to be applied for general model reduction purposes. To this aim, we pre-process all the data in this and the next section, by normalisation and sign splitting so that CDT and Radon-CDT are well defined. The normalisation constant is computed for all images such that the discrete sum is one. The output image is then rescaled appropriately before analysing the results. When interpolating, the normalisation constant for the unseen interpolated snapshot is also interpolated using the same interpolation used for the image.
Most of the examples presented here are strictly positive fields. However, when the initial image contains negative values, the images are split into positive and negative parts, and the transformation and subsequent steps (POD and interpolation) are applied on each part separately. The final image after inverse transformation is then obtained by summing up positive and negative reconstructions.

To test RCDT consistency, we present some initial test cases as one dimensional image-like inputs to the RCDT algorithm. These are transformed into RCDT space through the numerical algorithm, a discretisation of the continuous RCDT presented previously, before being converted back into the original `physical space'. Comparing the output image to the original input image then allows analysis of the numerical consistency of the RCDT algorithm, identifying any error introduced in the procedure and allowing us to observe the \textit{intrinsic} error caused by RCDT. Knowledge of this error is necessary for later studies involving model order reduction and predictive interpolation, so the \textit{intrinsic} error from the RCDT itself is already accounted for.

% %%%%%%%%%%%%%%%%%%%%%%%%%%%%%%
% %%%%%%%%%%%%%%%%%%%%%%%%%%%%%%
% \subsection{Numerical implementation}
% %%%%%%%%%%%%%%%%%%%%%%%%%%%%%%
% %%%%%%%%%%%%%%%%%%%%%%%%%%%%%%

% \SetKwComment{Comment}{/*}{*/}

We first analyse RCDTs' capability in retaining a single
snapshot image's geometrical shape, i.e. its consistency, observable by
successive forward and inverse RCDT on the given image. Comparing the
outputted forward and inverse image with the original, we focus on what we have defined as the \textit{intrinsic}
error -- a consequence of non-linear transformation -- that we desire. Finally,
the average relative Euclidean norm error, a discretised version of the relative $L^2$-norm error, is computed as a scalar quantifier to help
compare among image cases, defined as 
\begin{equation}\label{eqn:error}
    e_{2} = \frac{||x - \tilde{x}||}{||x||} = \frac{\left (\sum (x_i-\tilde{x_i})^2\right )^{1/2}}{\left (\sum(x_i^2) \right )^{1/2}},
\end{equation}
where $x_i$ is the true cell value and $\tilde{x_i}$ is the approximated cell value at cell $i$.

\subsection{Algorithm verification}\label{sec::RCDTverify}
\begin{figure}[htbp]
    \centering
    \includegraphics[width=\textwidth]{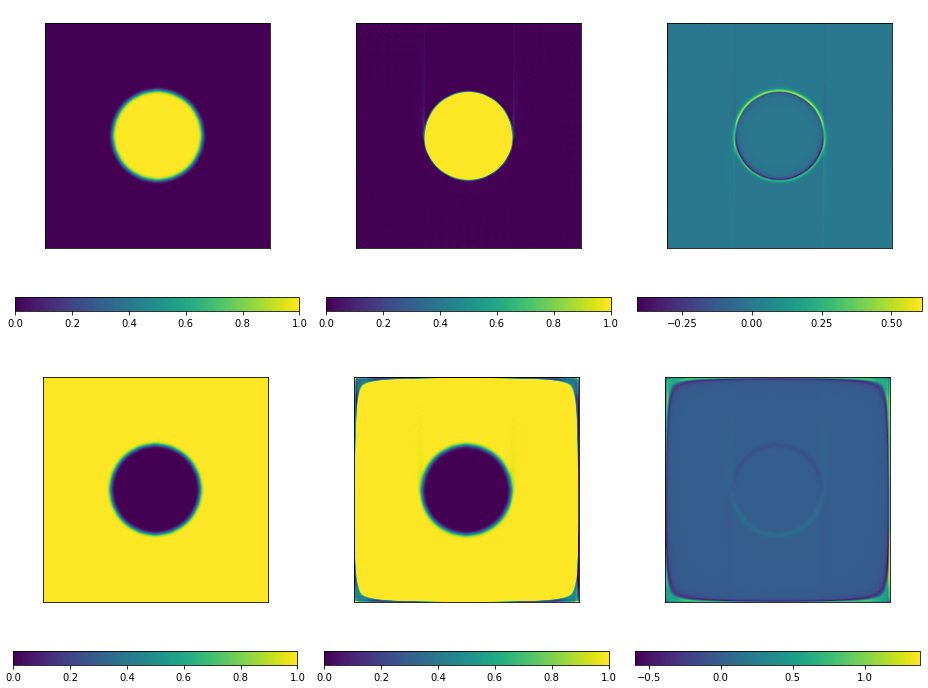}
    \caption{\small \textbf{Left}: input images of a smoothed circle defined by ones (top) or zeros (bottom) and vice versa for background. \textbf{Centre}:  result of the RCDT followed by iRCDT on the input images. \textbf{Right}: difference between the input images and RCDT iRCDT results. }
    \label{fig:circle_smoothed}
\end{figure}

\begin{figure}[htbp]
    \centering
    \includegraphics[width=\textwidth]{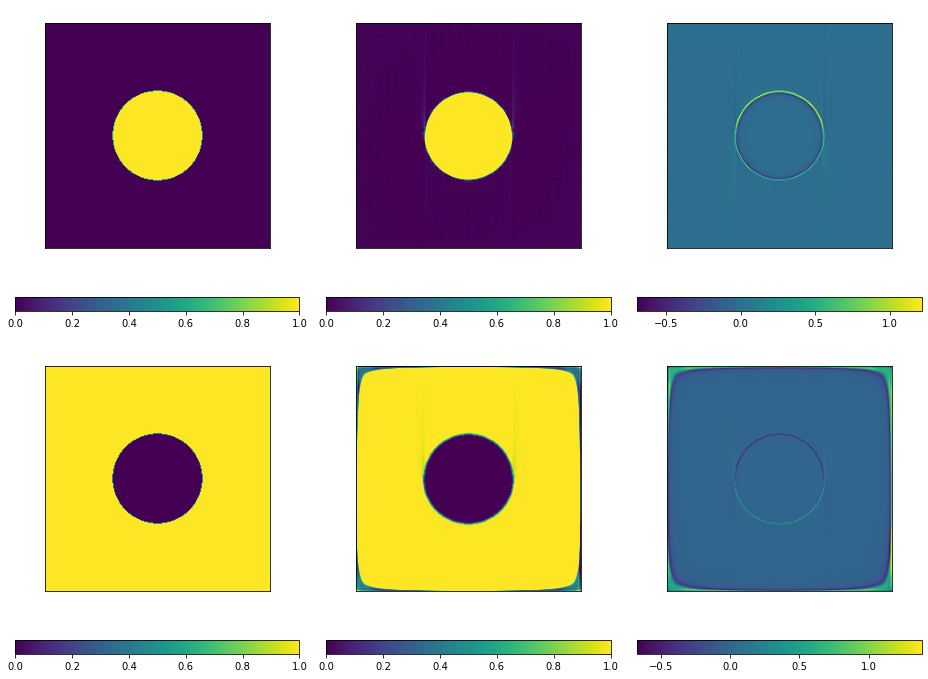}
    \caption{\small \textbf{Left:} input images of a circle defined by ones (top) or zeros (bottom) and vice versa for background. \textbf{Centre:} result of the RCDT followed by iRCDT on the input images. \textbf{Right:} difference between stated the input images and RCDT iRCDT results.}
    \label{fig:circle}
\end{figure}

\begin{figure}[htbp]
    \centering
    \includegraphics[width=\textwidth]{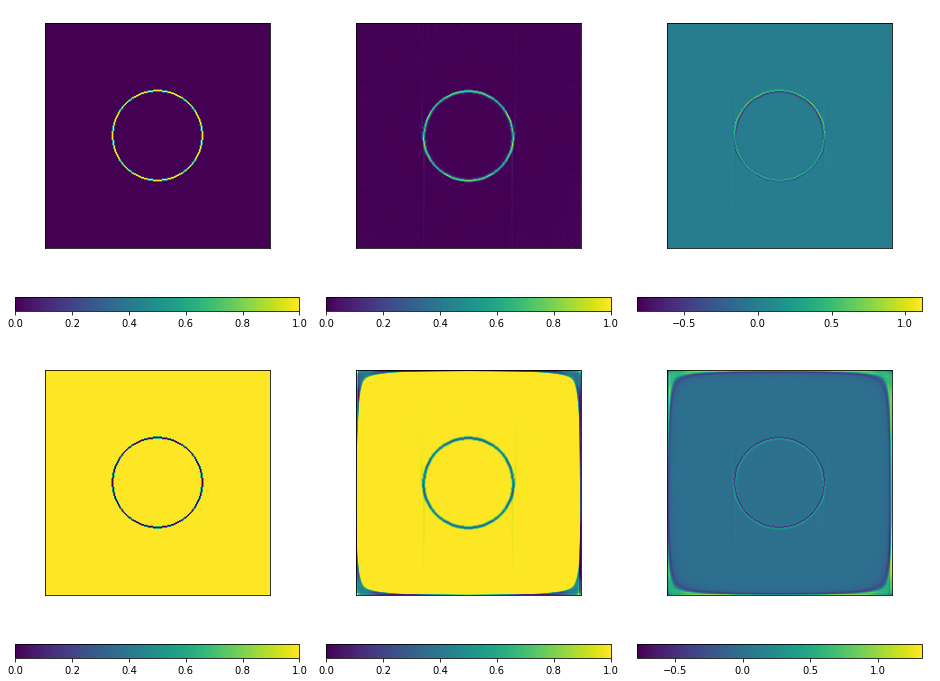}
    \caption{\small \textbf{Left:} input images of a circular edge ring defined by ones (top) or zeros (bottom) and vice versa for the background. \textbf{Centre:} result of RCDT followed by iRCDT on the input images. \textbf{Right:} difference between stated input images and RCDT iRCDT results.}
    \label{fig:edge}
\end{figure}

\begin{figure}[htbp]
    \centering
    \includegraphics[width=\textwidth]{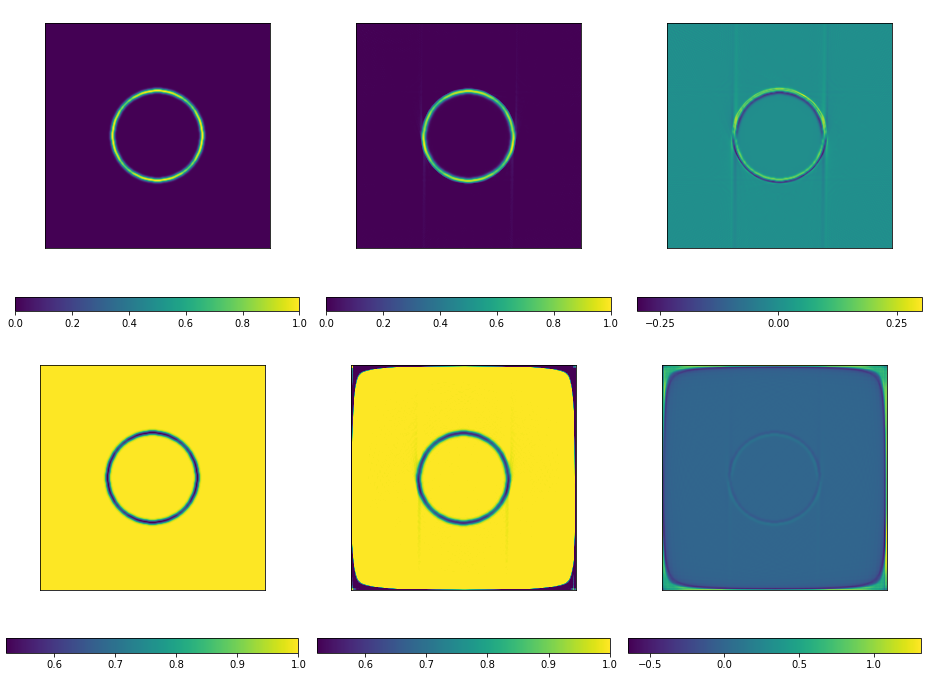}
    \caption{\small \textbf{Left}: input images of a smoothed circular edge ring defined by ones (top) or zeros (bottom) and vice versa for the background.. \textbf{Centre:} result of RCDT followed by iRCDT on the input images. \textbf{Right:} difference between stated input images and RCDT iRCDT results.}
    \label{fig:edge_smoothed}
\end{figure}

\begin{figure}[htbp]
    \centering
    \includegraphics[width=\textwidth]{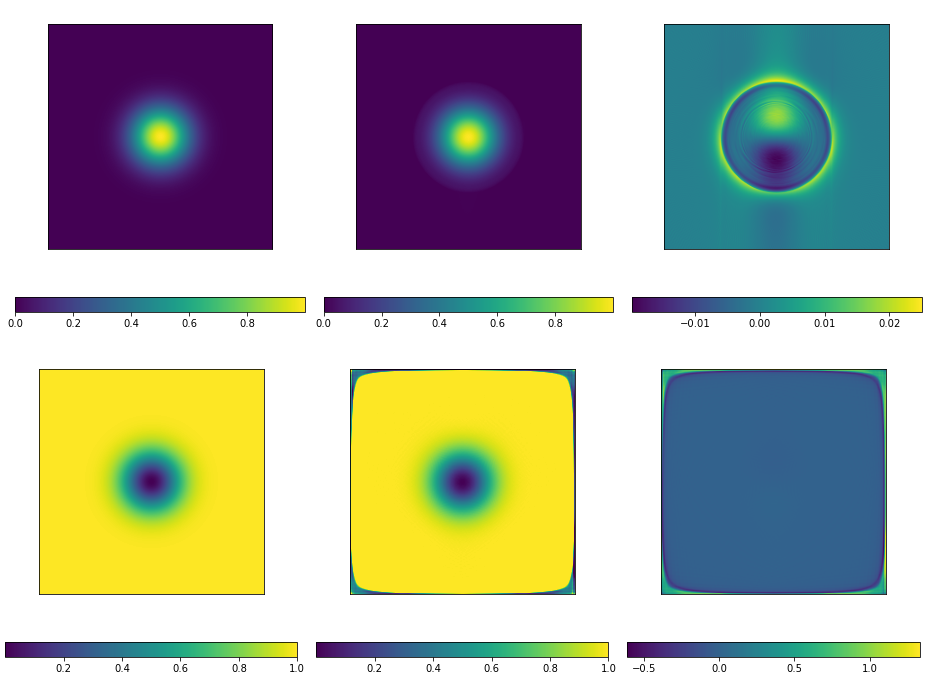}
    \caption{\small \textbf{Left:} input images of a normal Gaussian distribution function (top) ranging from 0 to 1, and vice versa (bottom). \textbf{Centre:} result of RCDT followed by iRCDT on the input images. \textbf{Right:} difference between stated input images and RCDT iRCDT results.}
    \label{fig:gaussian}
\end{figure}

%%%%%%%%%%%%%%%%%%%%%%%%%%%%%%
%%%%%%%%%%%%%%%%%%%%%%%%%%%%%%
Here, we present a selection of results focused on testing the consistency of the discrete RCDT algorithm and the \textit{intrinsic} error introduced. These test cases consist of simple images such as a smoothed circle with a radius of 50 pixels, \cref{fig:circle_smoothed}, its sharp counterpart,  \cref{fig:circle}, a ring obtained from applying an edge filter to a circle, \cref{fig:edge}, the smoothed counterpart case, \cref{fig:edge_smoothed}, and a Gaussian function, \cref{fig:gaussian}. The image values range between 0 and 1. Tests were also performed when these test images were `inverted'. That is to say, the `inverted' image pixels had values of $1 - x$ where $x$ is the original pixel value. In practice, this leads to the image features taking on a value of 0, whereas the background domain value is 1 (in the following referred to as the `inverted' images; for example, see \cref{fig:circle}). The RCDT was performed on the image for each test case, followed by the inverse RCDT (iRCDT) to return to the physical domain. The resulting image is then compared to the original input image; the difference between them is calculated cell-wise, allowing us to observe any \textit{intrinsic} error or, in other words, artefacts of the RCDT.

All test cases within this section of observing RCDT's \textit{intrinsic} error, i.e. image consistency, are performed on a uniform two-dimensional grid of $250$ by $250$ pixels (px); ergo $62500$~px in total. The various images used to test the \textit{intrinsic} error are constructed as $250$ by $250$ data arrays, each element representing a pixel and its 'intensity' -- e.g. flow velocity or an object -- within the array. For RCDT, the Radon transform projection angles and reference distribution of CDT are set to $180$ evenly spaced for all cases. Each test case shows some inconsistency in the RCDT, iRCDT procedure, as expected by a non-linear transform. In the outputted images, we observe an introduction of artefacts.

In the circle test case seen in \cref{fig:circle}, the resulting image differs from the input around the circle's edge. In particular, we observe an undershoot and overshoot pattern, with the inner edge of the circle showing a slightly lower predicted value and the outer edge showing higher errors. A similar effect is observed in the circle edge test image in \cref{fig:edge}; however, the extent of the error appears more localised, possibly due to the small width of the feature in the image.
In the inverted cases, another unique artefact is seen at the corners of the resulting images. For example, in \cref{fig:circle}, there is a relatively large overshooting of the background value in the corners of the image, extending partially along the edges of the image. This effect is visible in the inverted image cases and not in regular images.

In the Gaussian test case, \cref{fig:gaussian}, the error around the edges of the Gaussian feature is similar to the circle cases but with additional over and undershooting inside the feature. Despite this extra qualitative error, the quantitative error is an order of magnitude smaller than the circle test cases. This under/overshoot is possibly being caused by the change of values in the image domain, with a discontinuous and large value change being associated with a larger error; for example, the error values are much larger for the discontinuous circle test case compared to the continuous Gaussian test case. This is also supported by the smoothed circle and smoothed edge test cases, seen in \cref{fig:circle_smoothed,fig:edge_smoothed} respectively. In these cases, the discontinuous edges of the features are smoothed out with a Gaussian filter to render the changes in the image continuous. We see that both the under/overshoot error magnitudes and the overall image errors in these cases are smaller than their discontinuous counterparts, supporting the hypothesis that a smooth change in values is handled better by the transform and inverse transform algorithms and results in a more consistent result.

These errors arise from the numerical implementation of the CDT being used for the transformation rather than from the Radon transform algorithm. This was concluded after the images were transformed first into Radon space; this Radon space image was then iRadon-transformed, and this result was compared to the image after applying CDT-iCDT before iRadon once again. The Radon-iRadon image was free of artefacts, whereas the RCDT-iRCDT image contained significant edge artefacts. The result of the Radon space image undergoing CDT-iCDT back to Radon space was also compared with the original Radon-transformed image to see the errors introduced by the discrete CDT procedure in the Radon space.

\begin{table}
    \centering
    {\renewcommand{\arraystretch}{1.2}\begin{tabular}{ |l|c| }
        \hline
        \textbf{RCDT test cases} & \textbf{Relative} $\mathbf{L^2}$\textbf{-norm} \textbf{error}\\
        \hline\hline
        Circle & $2.279\times10^{-1}$  \\
        Circle inverse & $1.658\times10^{-1}$   \\
        Circle smoothed & $1.322\times10^{-1}$ \\
        Circle smoothed inverse  & $1.618\times10^{-1}$  \\
        Circle edge & $6.065\times10^{-1}$ \\
        Circle edge inverse & $1.708\times10^{-1}$  \\
        Circle edge smoothed & $2.356\times10^{-1}$  \\
        Circle edge smoothed inverse & $1.624\times10^{-1}$  \\
        Gaussian & $3.204\times10^{-2}$  \\
        Gaussian inverse & $1.628\times10^{-1}$  \\
        \hline
    \end{tabular}} 
    \caption{\small Relative $L^2$-norm error values for RCDT test cases. Quantitatively, errors for the inverted version of some cases fared better than their original counterparts, though marginally, contrasting to the qualitative differences seen in test case figures. Smoothed cases -- of the originals -- for the circle and circle edge also fared better by roughly factors of a half and third, respectively.}
    \label{tab:tests}
\end{table}

The error values for all RCDT test cases can be found in \cref{tab:tests}. Notably, the error for the standard Gaussian case is an order of magnitude lower than all other examples tested, whereas the least accurate reconstruction was in the circle edge case. Of our examples, the Gaussian case is the example where the inputs to the CDT function from the Radon transform algorithm most resemble a well-defined, smooth probability density function. Seeing as the CDT is defined from the space of smooth probability density functions, it is natural that this example returns the most consistent result. This is compared to the others where the input to the CDT has sharper features.

% From the features seen and errors encountered with these test examples, we can gain some insight into how helpful the RCDT may be when applied to CFD. Firstly, the fact that many examples of CFD data will not exhibit sharp discontinuities is beneficial for the RCDT because although it does deal with sharp boundaries, it is much less error-prone when used on smooth transitions and boundaries. Another feature seen is the rounding of edges and the addition of zero values along the non-zero borders of the inverted images. This could lead to problems when applied to CFD data and must be accounted for. This is only a problem when a large part of boundary values are non-zero, so this might not be an issue in some specific cases, e.g., interacting twin jet flows in \cref{fig:rcdt_prop_interp_test}.

%%%%%%%%%%%%%%%%%%%%%%%%%%%%%%
%%%%%%%%%%%%%%%%%%%%%%%%%%%%%%
\subsection{Interpolation study}\label{sec::prelimInterpErrStudy}
%%%%%%%%%%%%%%%%%%%%%%%%%%%%%%
%%%%%%%%%%%%%%%%%%%%%%%%%%%%%%

\begin{figure}[htbp]
    \centering
    \includegraphics[width=\textwidth]{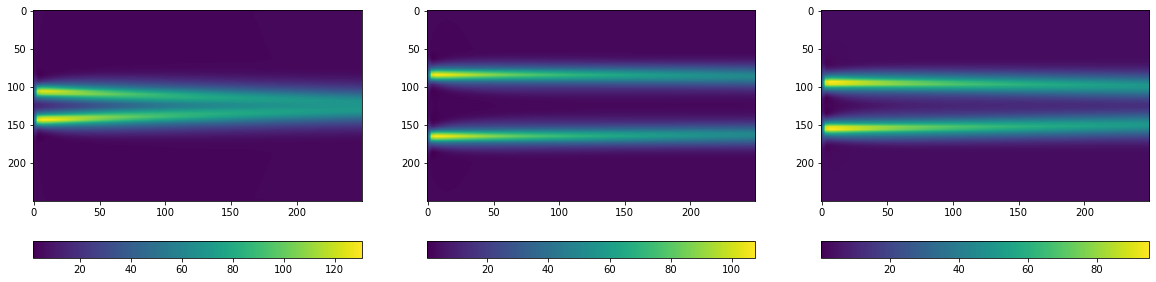}
    \caption{\small Input and target twin jet flow configuration images for differing separation widths. \textbf{Left}: Input image 1; Jets are close together; flows interact and merge close to sources. \textbf{Centre}: Input image 2; jet sources are further separated; therefore, flows do not merge together. \textbf{Right}: Target image; a jet configuration in between input images 1 \& 2.}
    \label{fig:rcdt_prop_inputs}
\end{figure}

\begin{figure}[htbp]
    \centering
    \includegraphics[width=\textwidth]{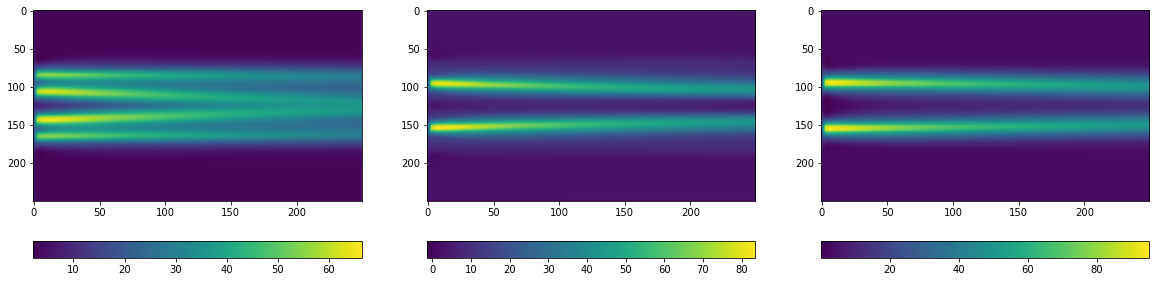}
    \caption{\small The resulting jet flow configuration images from interpolating the two input images in \cref{fig:rcdt_prop_inputs}. \textbf{Left}: Interpolation in physical space; failure of physical space interpolation. \textbf{Centre}: Interpolation in the RCDT space; qualitatively similar to target image, some over-shoot in inter-flow contact angle and smoothed jet boundaries. \textbf{Right}: Target image.}
    \label{fig:rcdt_prop_interp_test}
\end{figure}

% For studies of real-world applications, 'controlled' errors, like the \textit{intrinsic} error of RCDT above, can be a good compromise for the capability of interpolating highly nonlinear geometric flow features, a crucial aspect for modelling a variety of singular and multiple flow configurations for various parameter (spatial and transient) variations. This is partly due to how these 'controlled' errors can be compensated for -- or ultimately removed -- in post-processing. RCDT's \textit{intrinsic} error is mainly observed around the image and flow boundaries (see \cref{fig:circle,fig:edge} and later on \cref{fig:multi_ROM_20modes}) but still permits good image/flow capture.
%The same can be said of the 'controlled' POD reconstructive error mentioned in \cref{POD}; analysed later in \cref{Applications}. Considering a large number of spatial modes negates the error, but post-processing could again be utilised. Despite these caveats of RCDT -- and likewise ROM/POD -- in introducing systematical errors, the novel capability of interpolating nonlinear flow features makes it a small price to pay for such a valuable function.
The advantage of image interpolation in the RCDT space can be seen in the test case in \cref{fig:rcdt_prop_inputs,fig:rcdt_prop_interp_test}.
\Cref{fig:rcdt_prop_inputs} shows the flow field of two twin jet configurations, differing by the separation width of the jets and subsequent intertwining flows. These have been obtained by two-dimensional RANS turbulent flow simulations. The images are available in \cite{rcdt-rom-lib}.
The first two configurations, representing the largest and smallest separation widths, are interpolated in \cref{fig:rcdt_prop_interp_test} in real space (left) and RCDT space (middle) to obtain the target flow image on the right side.  The interpolation in real space results in a duplication of the jets. The RCDT one, instead, while having some over-shoot in inter-flow contact and dissipation, is significantly more accurate and realistic, preserving the number of jets and their structure.
Whilst only qualitatively shown for now, the interpolation error may over-shadow the \textit{intrinsic} one in transport-dominated problems, making it negligible in comparison.
%, similarly to the reconstruction error of POD/ROM, observed and tested later in \cref{Applications}; determined by the number of spatial modes taken into account in model reduction. The more modes, the better POD can reconstruct the image, so the less reconstruction error is seen in the output. It should be stated again that these 'controlled' errors can be compensated in image post-processing; later studies will utilise this to focus more on analysing the interpolation error briefly shown in \cref{fig:rcdt_prop_interp_test}, the end goal of accurately predicting flows using our RCDT-POD ROM workflow for a large number of parameter variations and flow configurations.

%%%%%%%%%%%%%%%%%%%%%%%%%%%%%%
%%%%%%%%%%%%%%%%%%%%%%%%%%%%%%
\section{RCDT-POD for model order reduction} \label{Applications}
%%%%%%%%%%%%%%%%%%%%%%%%%%%%%%
%%%%%%%%%%%%%%%%%%%%%%%%%%%%%%

In this section, we introduce the POD as a model reduction technique applied in RCDT space. The truncation performed by selecting a finite number of POD modes introduces a reconstruction error that is, in the $L_2$ sense, proportional to the singular values corresponding to the neglected modes. When POD modes are interpolated to obtain a ROM, the interpolation error also comes into play. The aim of this section is to study the effect of the RCDT on the reconstruction and interpolation errors and to compare the results with the standard POD in physical space.
% To understand this reconstruction error in the RCDT and physical space, we look towards comparison against high-fidelity CFD simulated data.
% The RCDT-POD procedure is compared against standard POD in physical space, and as such, the aforementioned \textit{intrinsic} error is present.
% To observe changes in the reconstruction error of POD, we compare with the original data for choices of $r=5$ modes and $r=20$ modes per subsequent case. In the numerical examples, POD is applied to the data in the RCDT space. This procedure requires defining a finite number of angles (in $2D$ or hyperplanes in $\mathbb{R}^n$), and the POD is applied to the discrete data. Apart from \cref{sec:GaussROM_RCDT}, an initial example where $r=5$ modes almost perfectly constructs the original image. In the following cases, we build the ROM using POD for the spatial mode decomposition and a Gaussian process regression (GPR) for the model prediction \cite{nguyen_peraire_2015, nguyen_peraire_2016}.

It should be noted that of the following examples, all but one, \cref{sec:GaussROM_RCDT}, use data from CFD simulations. This data on arbitrary, possibly non-uniform meshes can directly be used for physical space ROM. However, the RCDT software employed here requires a uniformly spaced grid (image-like data) for the input. This, therefore, requires the original CFD data to be uniformly interpolated across the domain. This was done using the \textit{scipy} interpolation function \textit{griddata}. Consequently, the RCDT-POD result is compared to the uniform grid input data, not the original  CFD data. 

%\subsection{Numerical implementation of RCDT-POD}

For the implementation of proper orthogonal decomposition (POD) and model order reduction (MOR) we use the \textit{EZyRB} package \cite{demo_ezyrb_2018}.
%See algorithm \ref{algo::RCDT&POD} for the algorithmic procedure for the RCDT-POD ROM workflow of given $K_i$ input images at snapshot times $t_i$ for $i=0,\dots, N_t$ and $N_t\geq 0$. For constructing the ROMs, we employ POD for snapshot decomposition using $N_r$ modes, given some input training data.
Different interpolation methods can be used in the POD space. In this work, we have limited our study to interpolation in one-dimensional parameter space, so the choice of interpolation schemes has been found to be not important, so we have used a linear interpolation as the default option. However, we have also tested other interpolation techniques in some test cases, such as the Radial Basis Function (RBF) method. In these cases, the default parameters used for a multiquadric kernel are: smoothness equal to zero (smoothness is meant as a paramater value greater than zero that increases the smoothness of
the approximation. 0 is for interpolation which means the function will
always go through the nodal points in this case), and shape parameter equal to one.

\subsection{Gaussian pulse}\label{sec:GaussROM_RCDT}
%%%%%%%%%%%%%%%%%%%%%%%%%%%%%%
%%%%%%%%%%%%%%%%%%%%%%%%%%%%%%

\begin{figure}[htbp]
    \centering
    \includegraphics[width=\textwidth]{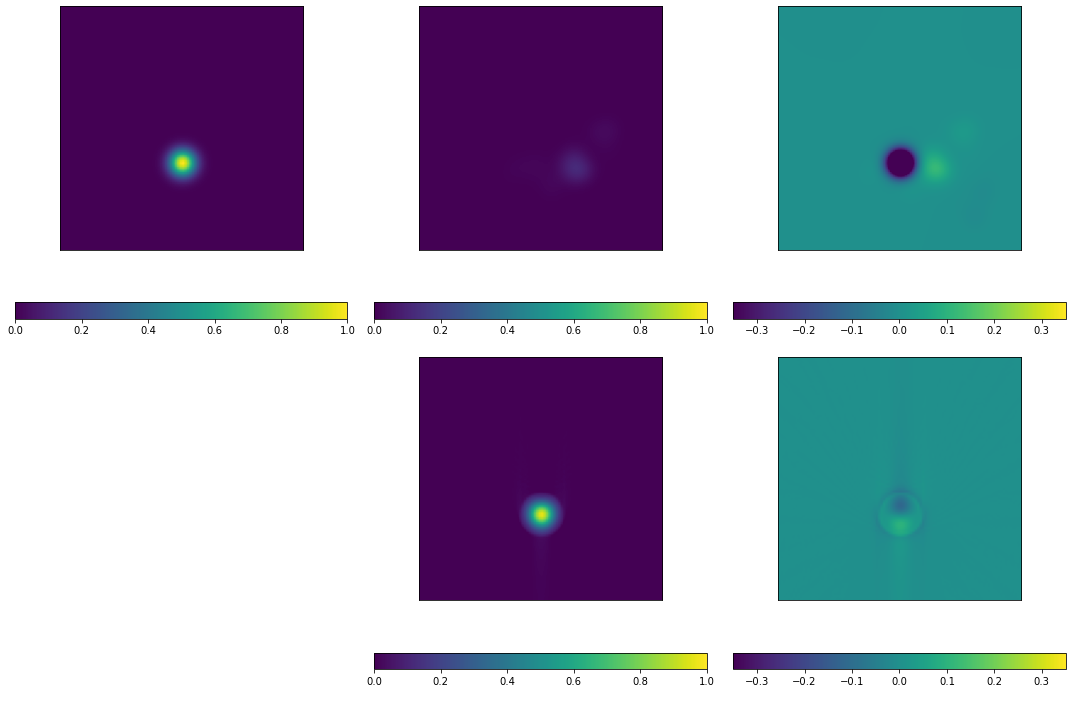}
    \caption{\small Reconstruction of the parameterised Gaussian pulse with POD (using the first five modes). \textbf{Top-left:} original snapshot image. \textbf{Top-middle:} physical space POD projection. \textbf{Bottom-middle:} RCDT-POD projection. \textbf{Top/bottom-right:} respective differences for each projection from original snapshot image.}
    \label{fig:Gauss_ROM}
\end{figure}

In \cref{fig:Gauss_ROM}, we study a Gaussian pulse in a 2D domain, similarly to what was investigated by \textit{Ren et al} \cite{ren_model_2021}
We consider a domain of size $100 \times 100$ cells and a Gaussian pulse parametrised by its position $(\mu_x,\mu_y$) in the $x$ and $y$ directions, and with a constant standard deviation $\sigma=5.3$. $\mu_{x,y}$ are randomly chosen in the interval $[10,90]$. $100$ snapshots have been used for the POD in the physical and the RCDT space. Considering only the first five POD modes, a random snapshot belonging to the training set is then reconstructed.

\begin{figure}[htbp]
    \centering
    \includegraphics[width=0.4\textwidth]{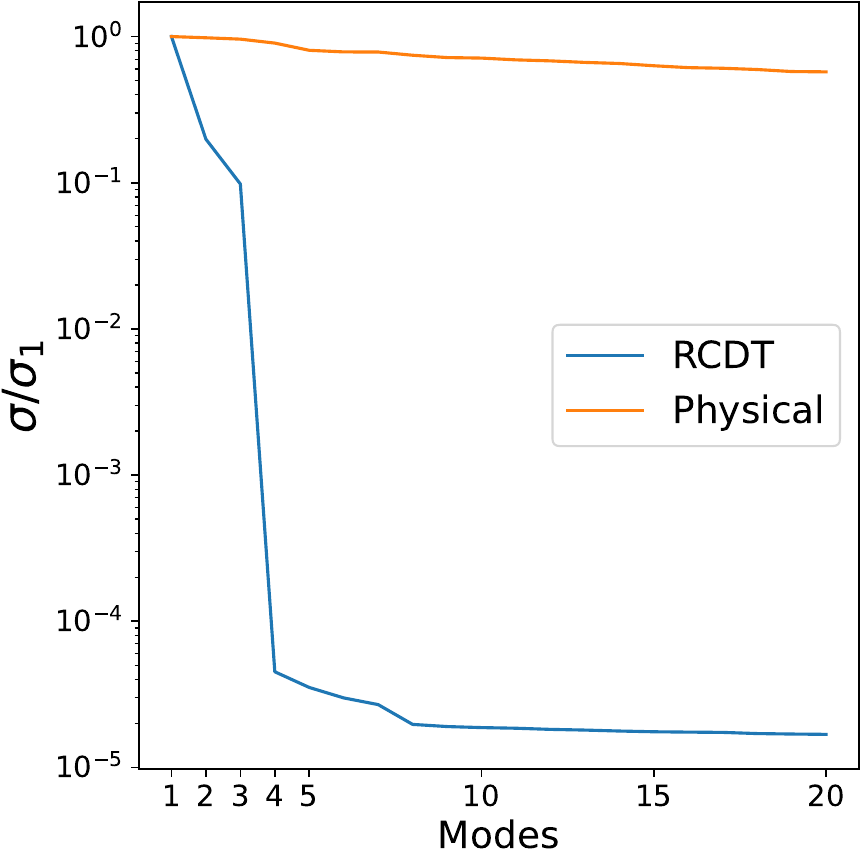}
    \caption{\small Plot of ratio of each singular value to first singular value in POD decomposition for physical space, and RCDT space POD's.}
    \label{fig:Gauss_ROM_Modes}
\end{figure}

In \cref{fig:Gauss_ROM_Modes}, we can see that the RCDT-based POD far exceeds the accuracy of the physical space-based POD, almost perfectly replicating the original snapshot. The number of POD modes required in RCDT space to capture most of the snapshots' energy accurately is considerably reduced with respect to the physical space case. This demonstrates the ability of the RCDT-POD procedure to compress travelling features efficiently. We can notice that, for each snapshot, we see an error pattern similar to the Gaussian cases presented previously, suggesting that most of the error here is due to the intrinsic RCDT error rather than the POD truncation.

In \cref{fig:Gauss_ROM_Modes}, the singular values of the POD decomposition (normalised by the first singular value) for the Gaussian POD are compared. Singular values represent the amount of `energy' stored in each corresponding POD mode (ordered from largest singular value to smallest). Therefore, a sharp decay in singular values corresponds to a larger proportion of the system's total energy captured by the first few POD modes, which, in turn, allows for a more accurate POD with fewer modes needed to capture the behaviour in the system. This example shows the physical space POD results in the first 20 POD mode singular values. This means that all of these POD modes contain a relatively large amount of the system's total energy, and hence, we would need to include many of the modes in the POD construction to get an accurate result. Conversely, in the case of POD in the RCDT space, we see a very sharp drop in singular values for the first four modes, as they contain the majority of the total energy of the system. As seen in \cref{fig:Gauss_ROM}, this results in an accurate POD using only the first five modes, whereas in physical space, the POD is inaccurate for the first five modes.

%%%%%%%%%%%%%%%%%%%%%%%%%%%%%%
%%%%%%%%%%%%%%%%%%%%%%%%%%%%%%
\subsection{Multi-phase wave}\label{sec:multiWaveROM_RCDT}
%%%%%%%%%%%%%%%%%%%%%%%%%%%%%%
%%%%%%%%%%%%%%%%%%%%%%%%%%%%%%

\begin{figure}[htbp]
    \centering    
    \begin{minipage}{0.49\textwidth}
        \centering       
        \includegraphics[width=\textwidth]{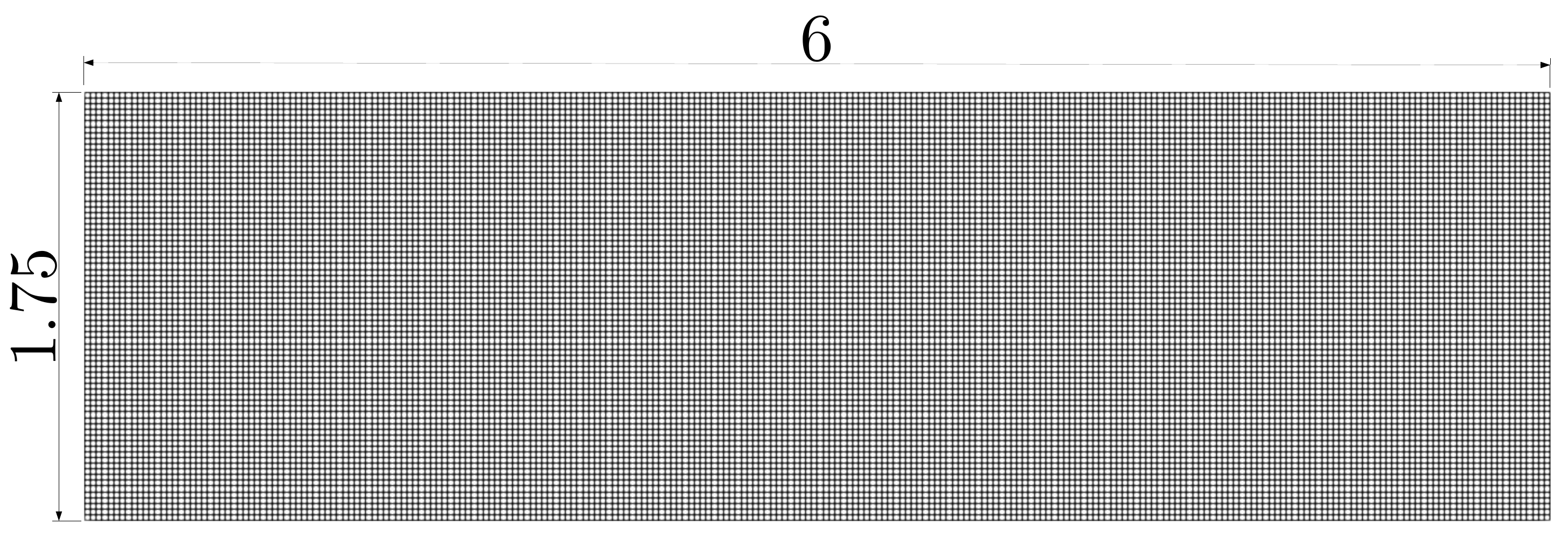}\\
        (a)  
    \end{minipage}      
    \begin{minipage}{0.49\textwidth}
        \centering    
        \includegraphics[width=\textwidth]{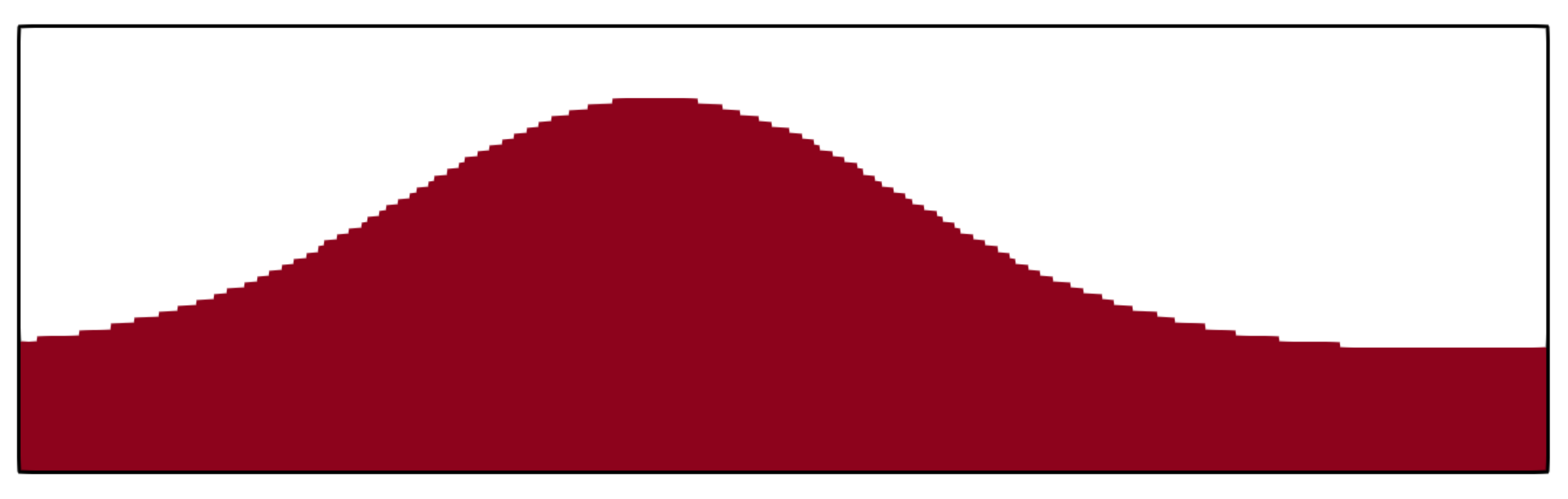}\\
        (b)    
    \end{minipage}
    \caption{Computational domain of the multi-phase CFD data set: (a) domain size and mesh grid structure, (b) initial condition for the $\alpha_w$ field}
    \label{fig:multi-phase_domain}
\end{figure}

The data set for this example has been generated by solving the unsteady Navier Stokes equations for two incompressible, isothermal immiscible fluids. The domain is given by the rectangle $\Omega = [-2.5,3.5]\times[-0.5,1.2]$ and the computational domain (see~\cref{fig:multi-phase_domain}) is obtained with a structured grid composed by $250 \times 75$ quadrilateral cells. The results are obtained using the \texttt{interFoam} solver developed in the OpenFOAM\textsuperscript{\textregistered} finite volume library. We refer to \cite{Larsen2019} for more details on the specific equations that are solved and for the numerical implementation.  For what concerns the velocity field, we apply a uniform velocity $U=(0.25,0) \unit{m/s}$ on the left side of the domain, a \texttt{noslip} condition on the bottom of the domain, and a \texttt{zeroGradient} (i.e. $\nabla U \cdot n = 0)$ condition on the top and right side of the domain. For the pressure field, we used a uniform \texttt{totalPressure} with $p=0 \unit{kg/(m s^2)}$ to the bottom boundary and a uniform \texttt{fixedFluxPressure} with $p=0 \unit{kg/(m s^2)}$ on other boundaries. The field $\alpha_w$, which varies between $0$ and $1$ and represents the fraction of volume of water in each cell, has a uniform condition $\alpha_w=0$ on the bottom boundary and a \texttt{zeroGradient} condition on the other boundaries. The properties of the two phases are set as $\rho_1 = 10^3 \unit{kg/m^3}$, $\nu_1 = 10^{-6} \unit{m^2/s}$, and $\rho_2 = 1\unit{kg/m^3}$, $\nu_2 = 1.48 \cdot 10^{-5} \unit{m^2/s}$, where $\rho_1$, $\rho_2$ and $\nu_1$, $\nu_2$ are the density and the kinematic viscosity of the first and second phase, respectively. Snapshots are collected in the $[0,5]\unit{s}$ time window.

For RCDT-POD, the CFD data is interpolated onto a uniform grid of size $200 \times 150$. The number of time steps in this data set is $N_t = 199$, discounting the initial. The initial condition for the velocity field is uniform $U=(0.25,0) \unit{m/s}$, for the pressure field is uniform $p=0\unit{kg/(m s^2)}$, while for the $\alpha_w$ field is set according to the following profile:

\begin{equation}
    \begin{cases}
        \alpha_w = 1 \quad \text{if} \quad y<e^{-0.5x^2},\\
        \alpha_w = 0 \quad \text{otherwise}.
    \end{cases}    
\end{equation}
The time integration is performed using an implicit first-order Euler scheme. The data set exhibits the wave travelling across the domain from left to right during the time series.

\begin{figure}[htbp]
    \centering
    \includegraphics[width=\textwidth]{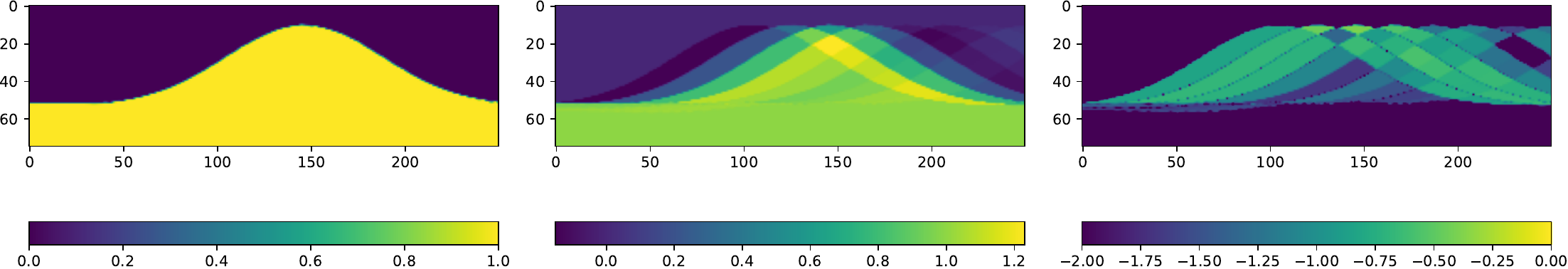}\\
    \includegraphics[width=\textwidth]{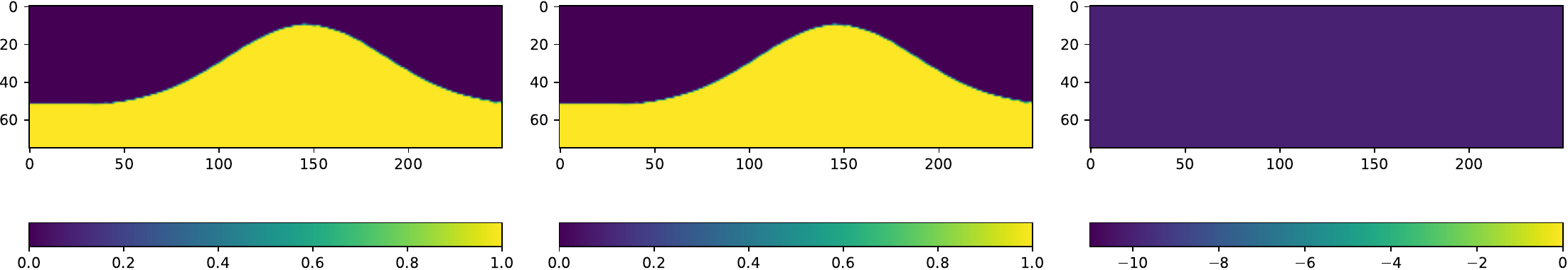}\\
    \includegraphics[width=\textwidth]{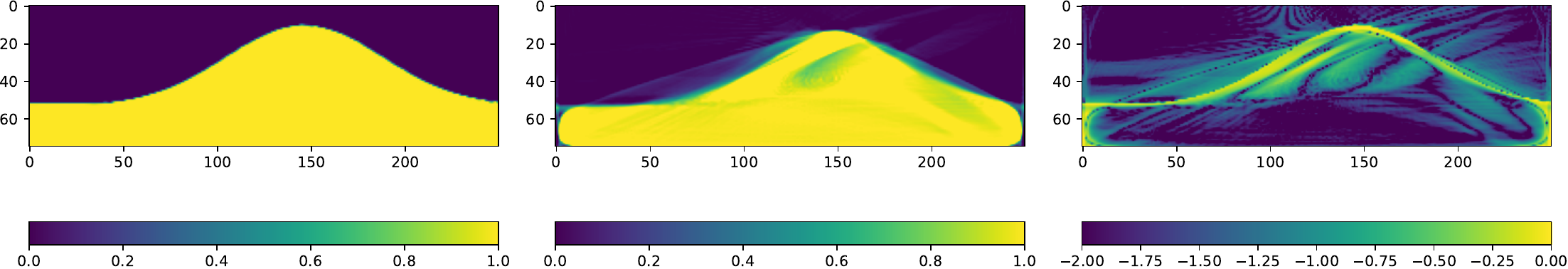}\\
    \includegraphics[width=\textwidth]{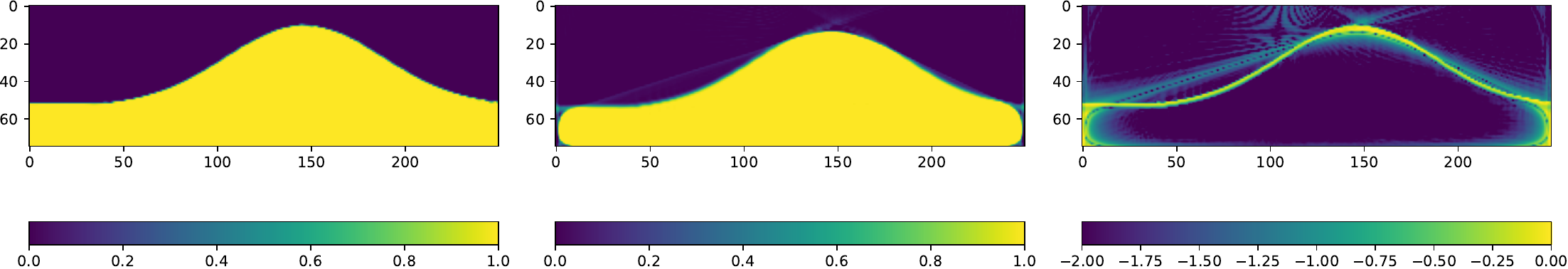}
    \caption{\small  Reconstruction of the multi-phase wave POD. Left: original image. Middle: reconstruction. Right: logarithm of the absolute error. From top to bottom, respectively, physical space POD with five modes, physical space POD with ten modes, RCDT-POD with five modes, and RCDT-POD with ten modes.}
    \label{fig:multi_ROM}
\end{figure}

\begin{figure}[htbp]
    \centering
    \includegraphics[width=0.4\textwidth]{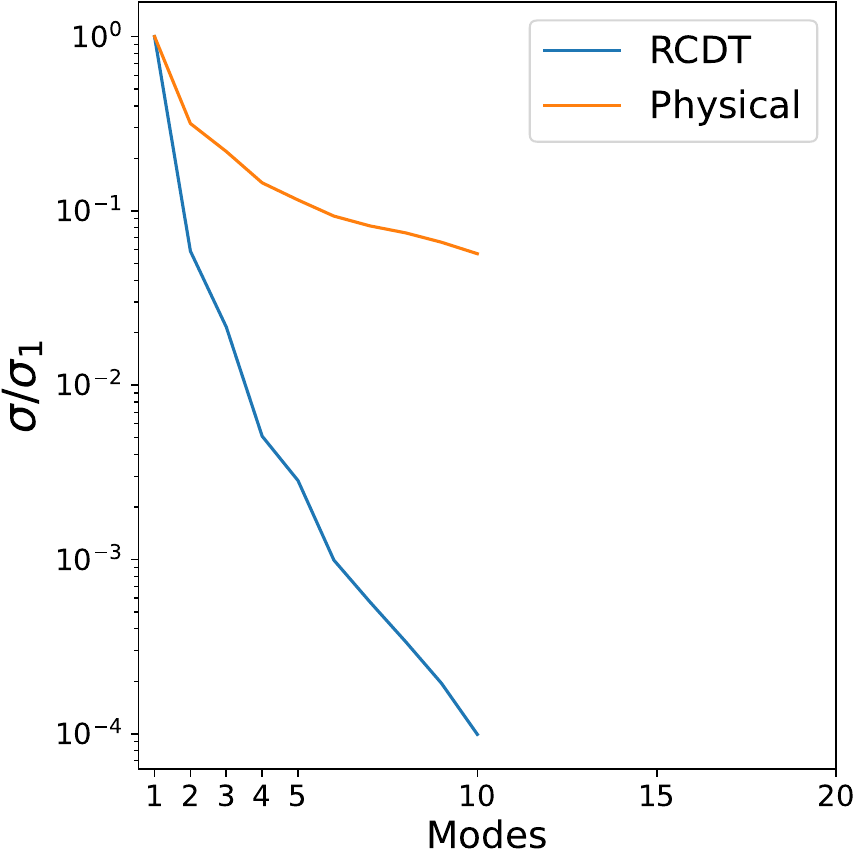}
    \caption{\small Multi-phase wave test case: singular values' decay for physical, and RCDT space.}
    \label{fig:multi_ROM_Modes}
\end{figure}

\Cref{fig:multi_ROM} shows the reconstruction of a single snapshot of the volume fraction in physical space (first two rows) and RCDT space (third and fourth row) with five (first and third row), and ten (second and fourth row) POD modes. The original snapshot is shown on the left, the reconstructed image in the centre, and the logarithm of the absolute error in the third column. There is a considerable difference in the physical space POD between five and ten modes. With ten modes, as expected, the physical space POD perfectly reproduces the image; however, when reducing the number of modes, the physical space POD struggles to capture the shape of the original wave. This is in contrast to the RCDT-POD, which retains more of the shape of the original when only five modes are used for reconstruction, although affected by the intrinsic RCDT error. The result improves slightly for ten modes, but most of the shape is captured in the first few modes, unlike in the case of physical space. This is made clear in \cref{fig:multi_ROM_Modes} where the decay of the POD singular values for the RCDT-POD case is much quicker than in the physical space case. However, this does not necessarily result in a lower overall error as it does not consider the intrinsic error in the RCDT/iRCDT transform. 

% The singular values' decay measures how well (in terms of $L^2$ error) we can approximate elements of the discrete solution manifold used for the training stage (the snapshots matrix) with a linear subspace of the prescribed dimension. It does now tell how well we can approximate the real solution manifold. 

The RCDT-POD results, in fact, see the appearance of edge-bound errors around the wave along the boundaries of the domain, especially in the bottom corners. Future work on this example will test the effect of adding a padding region around the data or using other forms of pre- and post-processing, such as thresholds and edge detection, which may improve results in certain cases. In conclusion, RCDT is usually better at capturing the geometric features and reducing the number of modes required to capture those. However, it might introduce non-physical geometrical features. We suppose that this fact might be caused by the intrinsic errors. More numerical experiments should be performed to have a clearer picture of the method's positive and negative features.

%%%%%%%%%%%%%%%%%%%%%%%%%%%%%%
%%%%%%%%%%%%%%%%%%%%%%%%%%%%%%
\subsubsection{Predictive time interpolation}\label{sec::multiWaveROM_RCDTPredict}
%%%%%%%%%%%%%%%%%%%%%%%%%%%%%%
%%%%%%%%%%%%%%%%%%%%%%%%%%%%%%

% Summary of the sections motivation
Having checked the RCDT-POD reconstruction on single snapshots, we study here the full RCDT-POD ROM workflow for predictive interpolation in time.  The case settings we consider are similar to those in \cref{sec:multiWaveROM_RCDT}. The difference here is that before inverting results from RCDT-POD space into the physical space, snapshots determined by the RCDT-POD workflow are linearly interpolated in time whilst still in RCDT-POD space. 
% Case setting of multi-phase case

\begin{figure}[htbp]
    \centering
    % \begin{subfigure}{\textwidth}
    %     \centering
    \includegraphics[width=\linewidth]{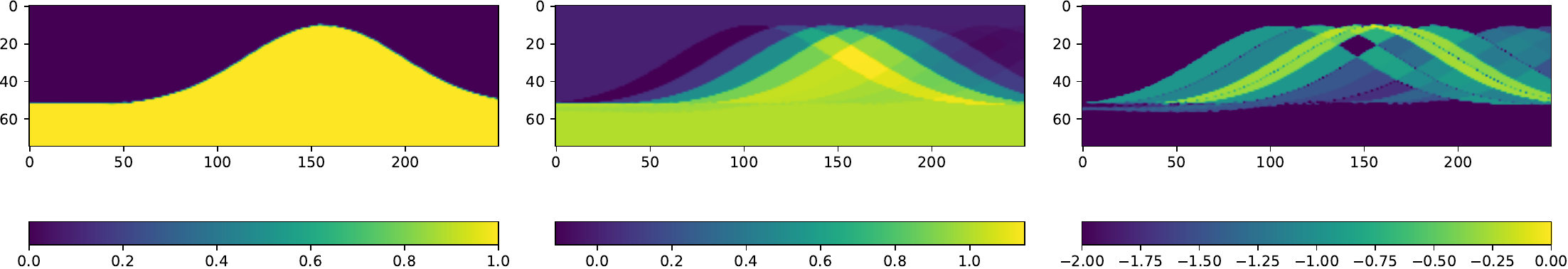}\\
    \includegraphics[width=\linewidth]{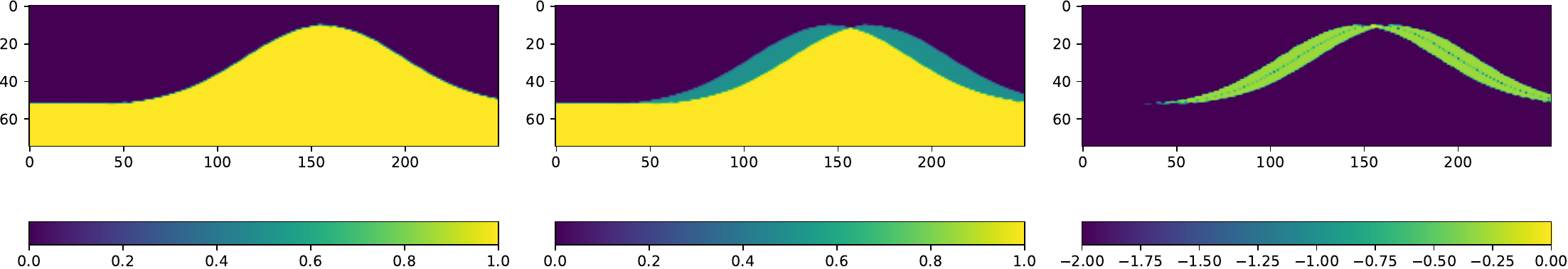}\\
    \includegraphics[width=\linewidth]{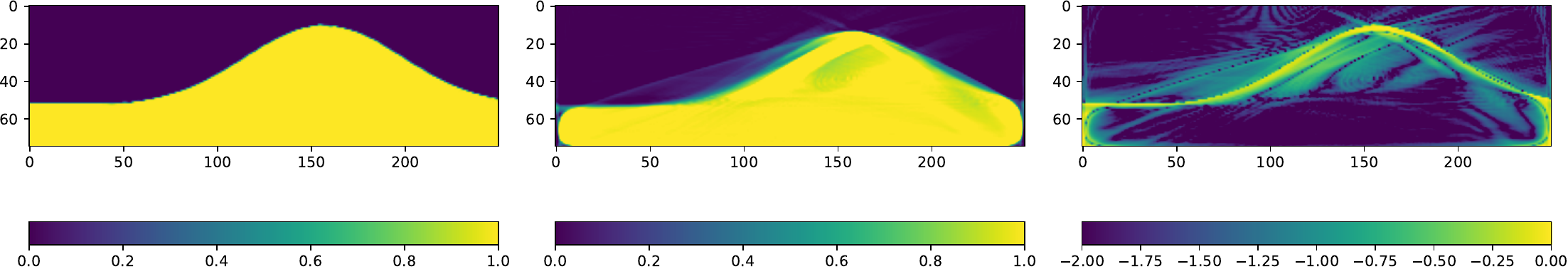}\\
    \includegraphics[width=\linewidth]{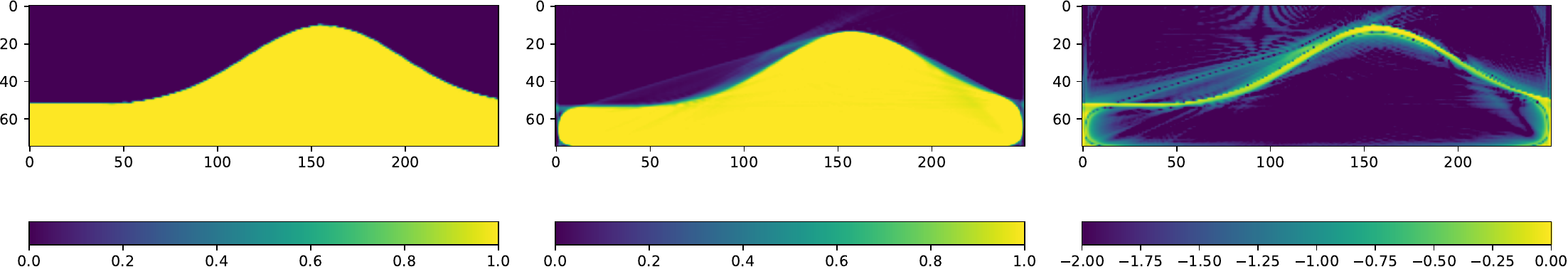}
        \caption{Interpolation of an unseen image of the multi-phase wave with POD. Left: original image. Middle: reconstruction. Right: logarithm of the absolute error. From top to bottom, respectively, physical space POD with five modes, physical space POD with ten modes, RCDT-POD with five modes, and RCDT-POD with ten modes.}
    \label{fig:interpTestMultiphase_Master}
\end{figure}

All transformations, decompositions, and interpolations are done using a reduced training data set of the original multi-phase CFD data. This reduced data $D_r$ set is composed by $10$ snapshots selected uniformly in the original $N_t=199$ snapshots $D_f$ (i.e. $D_r = D_f(1:20:\text{end})$ using matlab notation). The data is given on a $250$ by $75$ grid resolution. Linear interpolation is performed in each space at the chosen target snapshot index $k=75$ so as not to be a part of our training data set and compared against the original $k=75$ snapshot multi-phase data.
%Since the original data is made of 0-1 values, results are run through a threshold of 0.5, i.e., values above 0.5 are set to 1; otherwise, they are set to 0. This allows us to notice better the accuracy in detecting the wave shape.

Within \cref{fig:interpTestMultiphase_Master} are the results of the RCDT-ROM (i.e., interpolation in the RCDT-POD space) for five and ten modes, respectively. The left column shows the target snapshot; the middle column shows the predicted snapshot and the right column shows the logarithm of the absolute error between the target and predicted snapshots. The results of the physical space POD are shown in the top two rows, and the RCDT-POD results are shown in the bottom two rows. The results show that the RCDT-POD ROM can, from a qualitative standpoint, predict the target snapshot, even with only five modes. The physical space POD, however, is unable to accurately predict the target snapshot with only five modes, and even with ten modes, the prediction is not as accurate as the RCDT-POD ROM (for what concerns the wave shape). This is due to the RCDT-POD ROM being able to capture the wave shape with fewer modes than the physical space POD while
results in physical space show a failure of interpolation in the physical space, forming a bimodal wave with an inaccurate peak location to the target snapshot. However, we highlight that the RCDT-POD ROM introduces spurious conceptional boundary treatment problems and introduces unphysical geometrical features. Stanard POD is good in these respects.
% fftROM
% Within \cref{fig::fftROMTest} are the results of following the same procedure as in \cref{fig::physROMtest}, but performed in Fourier space instead. Like before, the figures show a distinct failure in interpolating the multi-phase wave within Fourier space. Unlike \cref{fig::physROMtest}, we observe clustering of artefacts on wave boundary, partly due to the thresholding of results and Fouriers' wave-like nature.
% rcdtROM
% Inside \cref{fig::rcdtROMTest} are interpolation results within RCDT-POD space. As can be seen by the respective predicted snapshot and scaled absolute error image, compared to \cref{fig::physROMtest,fig::fftROMTest}, there is a notable qualitative difference.
The RCDT prediction gives a much more accurate location of the wave peak, and the error at the wave boundary is primarily situated at the wave tail end. A consequence, however, as seen in all other tests of RCDT, is some rounding of the outer boundaries is present.
% rcdtInvert

% In \cref{fig::rcdtInverseTest} we provide a verification of RCDT's accuracy when transforming, and subsequently inversely transforming, the original input snapshot data back and forth from RCDT space. The target snapshot in \cref{fig::rcdtInverseTest} is the $50^{\text{th}}$ snapshot. Error is observed along the wave boundaries, including the external 'boundary' where we observe the rounding much like in \cref{fig::rcdtROMTest}.

%%%%%%%%%%%%%%%%%%%%%%%%%%%%%%
%%%%%%%%%%%%%%%%%%%%%%%%%%%%%%
\subsection{Airfoil CFD data set}\label{sec:airfoil}
%%%%%%%%%%%%%%%%%%%%%%%%%%%%%%
%%%%%%%%%%%%%%%%%%%%%%%%%%%%%%

\begin{figure}[htbp]
    \centering    
    \begin{minipage}{0.49\textwidth}
        \centering       
        \includegraphics[width=\textwidth]{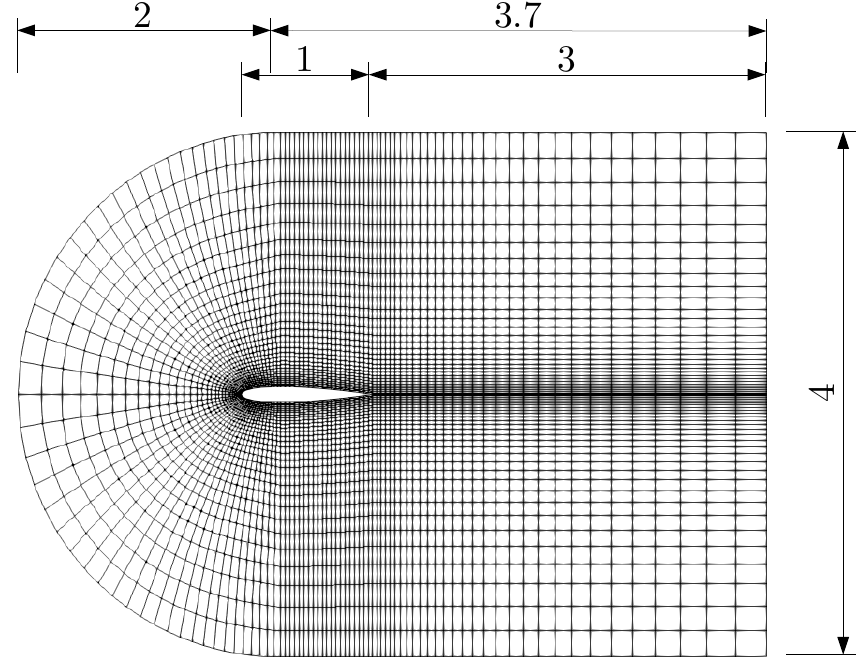}\\
        (a)  
    \end{minipage}      
    \begin{minipage}{0.49\textwidth}
        \centering    
        \includegraphics[width=\textwidth]{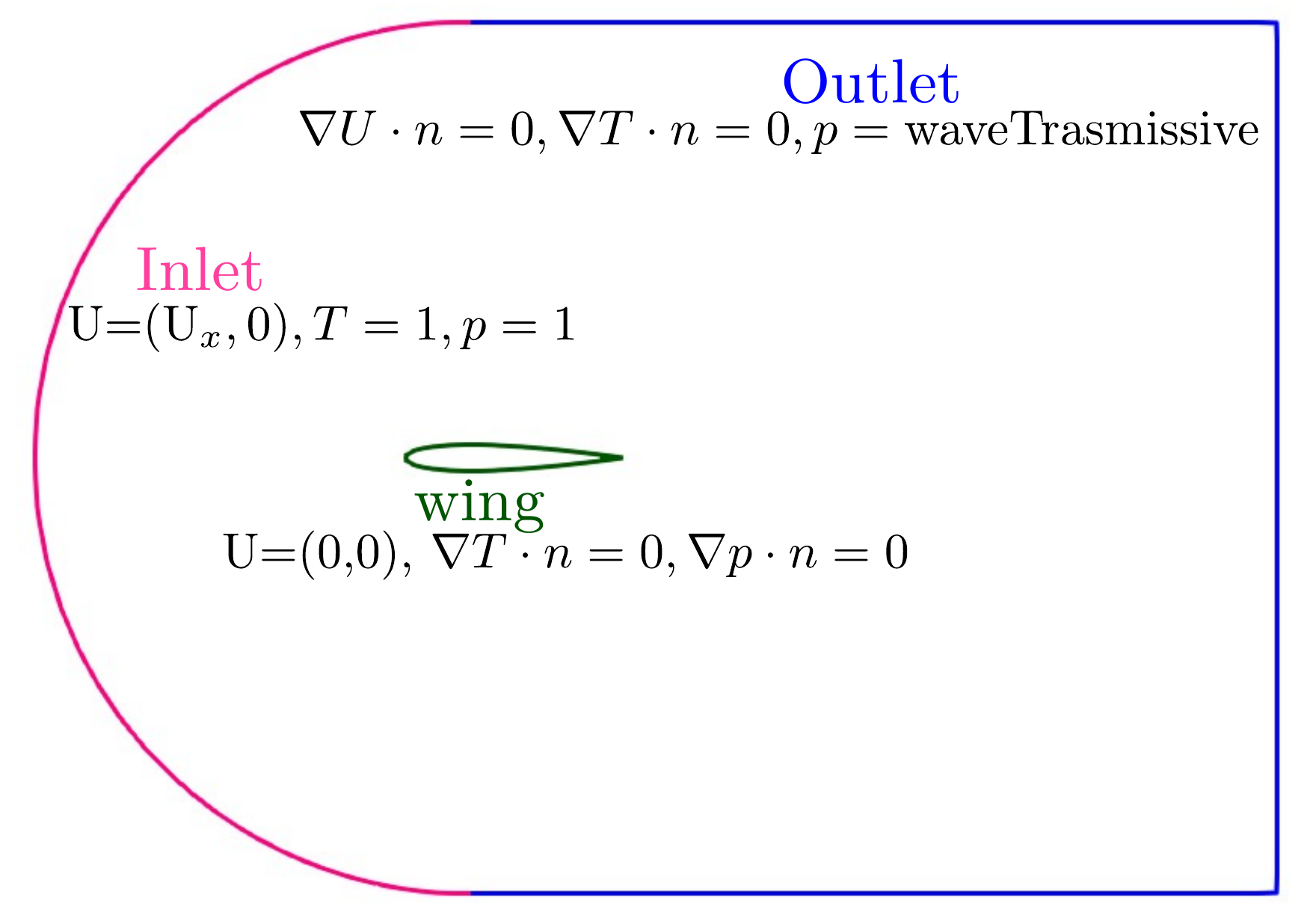}\\
        (b)    
    \end{minipage}
    \caption{Computational domain of the airfoil CFD data set: (a) domain size and mesh grid structure, (b) domain patches and boundary conditions.}
    \label{fig:airfoil_domain}
\end{figure}

In this test case, the data set is created by solving the compressible Navier-Stokes equations around a NACA0012  profile. The mesh is structured and counts $4500$ quadrilateral elements. For more details on the mesh structure and domain dimension, see \cref{fig:airfoil_domain}. The equations are solved employing a cell-centred finite volume method and using the \texttt{sonicFoam} solver \cite{Nakao2014} developed in the OpenFOAM library \cite{jasak1996error}. 

% original data
The simulation is 2D, and boundary conditions are set according to~\cref{fig:airfoil_domain}. In particular, the \texttt{waveTrasmissive} boundary condition is a special non-reflecting condition implemented in OpenFOAM. In the picture, $U$, $T$ and $p$ denote the velocity, temperature and pressure fields, respectively. The problem is transient; the time integration is carried out with an implicit Euler scheme, and the simulation time window is the interval $[0,3]$s. The time step used to solve the Navier-Stokes equations numerically is $\Delta t=0.001 \unit{s}$, and we store every ten time-steps. This means we acquire $300$ snapshots to test our numerical pipeline. Each snapshot, i.e. state, contains only the solved velocity profile to test our methodology. The inlet velocity $U_x$ is indirectly parameterised through the Mach number $\text{Ma}$, which is varied uniformly from $1.9$ to $3.1$ with intervals of $0.1$.

% \begin{figure}[htbp]
%     \centering
%     \includegraphics[width=\textwidth]{images/Airfoil_intrinsic.pdf}
%     \caption{\small Reconstruction of the airfoil magnitude velocity field (using the first five  modes). Left: original snapshot image. Middle: RCDT-POD space reconstruction. Right: logarithm of the absolute error.}
%     \label{fig:airfoil_intrinsic}
% \end{figure}

\begin{figure}[htbp]
    \centering
    \includegraphics[width=\textwidth]{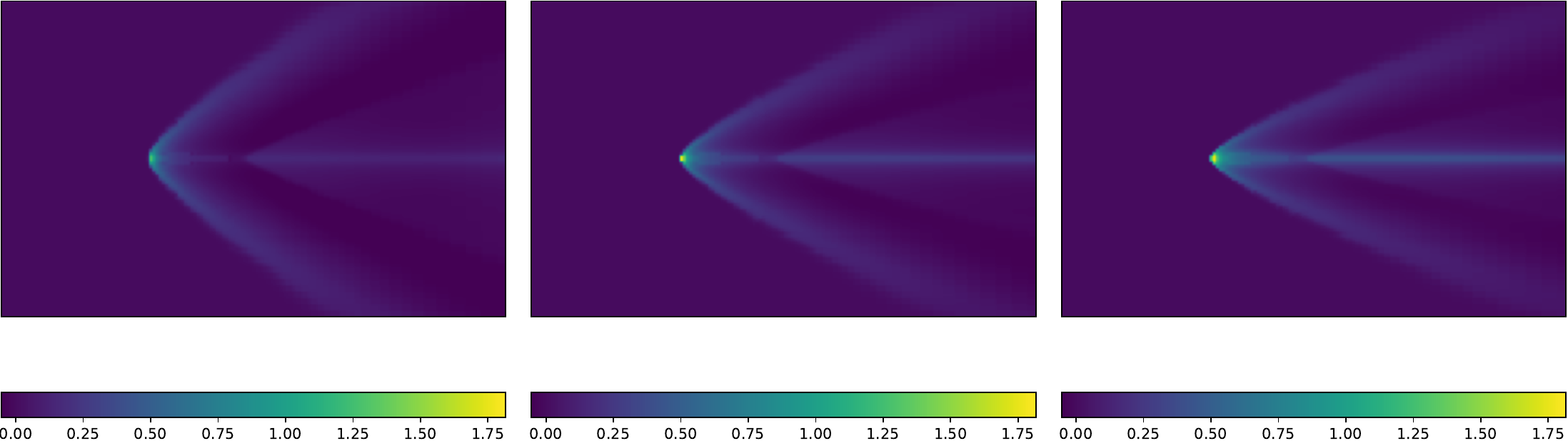}
    \caption{\small Airfoil test case. Right/Left: snapshots used. Middle: target snapshot.}
    \label{fig:airfoil_snapshots}
\end{figure}

\begin{figure}[htbp]
    \centering
    \includegraphics[width=\textwidth]{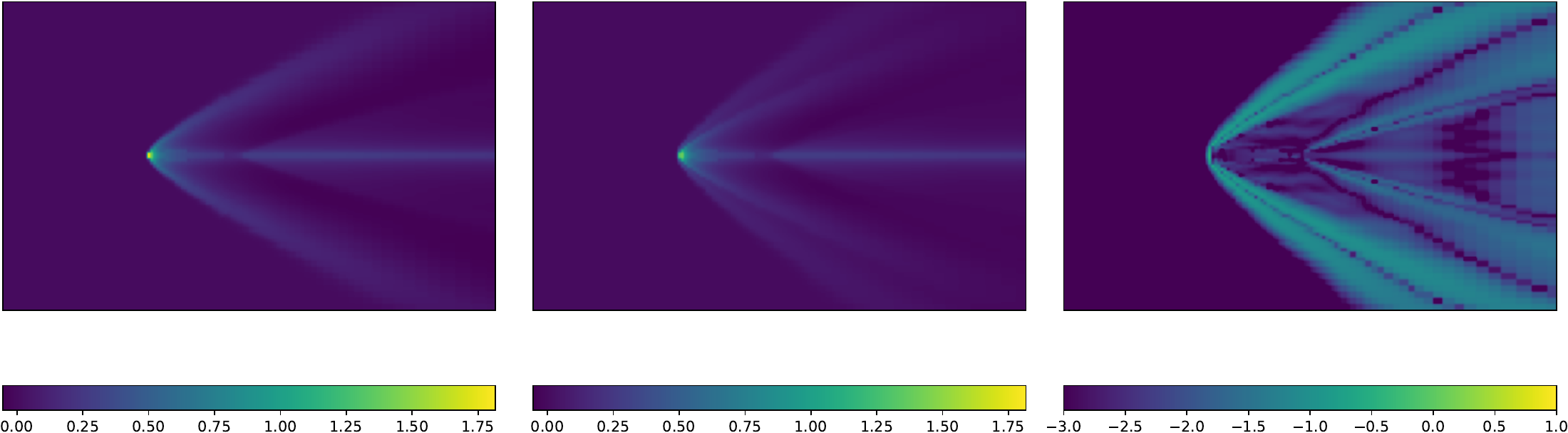}
    \includegraphics[width=\textwidth]{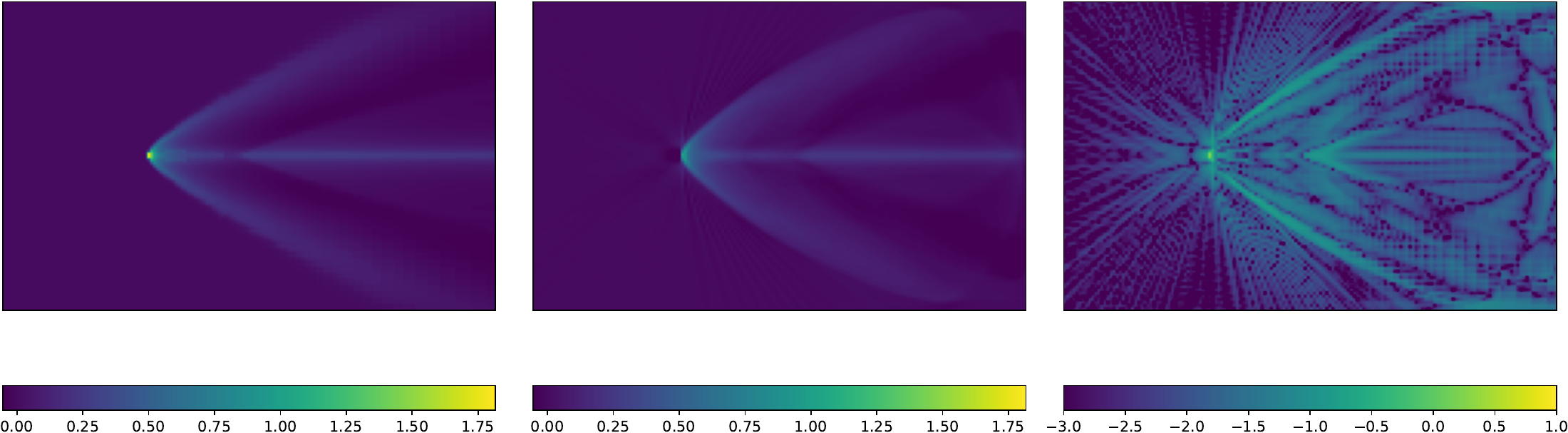}
\caption{\small Interpolation of an unseen image of the airfoil magnitude velocity field. Left: original image. Middle: reconstruction. Right: error. Physical space (top) and RCDT-POD space (bottom).}
    \label{fig:airfoil_ROM}
\end{figure}

% preprocessing
Here, we apply the ROM procedure to the relative velocity magnitude obtained by subtracting the inlet boundary velocity from the original velocity field. This pre-processing helps to avoid issues in performing the RCDT on fields that are not zero on the boundary. The CFD data is interpolated onto a uniform grid of size $250 \times 200$. The values on the regular rectangular grid are obtained via interpolation with the nearest point approach for what concerns extrapolation. This means that extrapolation in the area inside the rectangular grid but outside the computational domain is performed by assigning the value equal to the one at the nearest point inside the computational domain.
We use a signed variant of RCDT, applying the RCDT-POD to the positive and negative parts of the inputted CFD data separately and later re-combined. See \cite{aldroubi_signed_2022} for further details on the transformation of signed measures in RCDT space. This results in a double number of POD modes being used in interpolation and reconstruction.

% training snapshots
We present results on the interpolation in the Mach number space, considering only the last time-step. The chosen target snapshot is at $\Ma = 2.5$, and the training images are taken at $\Ma = 1.9$ and $\Ma = 3.1$.
In this test case, POD does not play a role, as the number of snapshots is equal to the number of modes used.
% \Cref{fig:airfoil_intrinsic} shows the intrinsic and reconstruction error of a single snapshot. This is magnified by the logarithmic scale.
\Cref{fig:airfoil_snapshots} shows the snapshots used subsequently for the interpolation. The left and right are the training snapshots for $\Ma=1.9$ and $3.1$, respectively, while the centre is our target snapshot for $\Ma=2.5$, which is not part of the training. The two snapshots used for the training are shown on the left and right.
The interpolation results are shown in \cref{fig:airfoil_ROM}. The physical space results are shown in the top row, and the RCDT-POD results are shown in the bottom row. The results show that the RCDT-POD ROM, despite the intrinsic error, can predict the target snapshot,  without introducing an additional shock within the wake, clearly visible in the physical space interpolation. However, unphysical and uncausal features upstream the airfoil contact are introduced by the RCDT-POD ROM. The classical POD has error zero in those regions and is physically more plausible.

\section{Conclusions} \label{Conclusion}
%%%%%%%%%%%%%%%%%%%%%%%%%%%%%%
%%%%%%%%%%%%%%%%%%%%%%%%%%%%%%
%%%%%%%%%%%%%%%%%%%%%%%%%%%%%%
%%%%%%%%%%%%%%%%%%%%%%%%%%%%%%

This work has focused on implementing and verifying the Radon-Cumulative Distribution Transform (RCDT) for image and flow capture and assessing its applicability in model order reduction (MOR) -- under proper orthogonal decomposition (POD) -- of high-fidelity CFD input data. RCDT and subsequent RCDT-POD MOR workflows were tested for accuracy compared against either the original input images or standard POD in physical space.
Although the results presented here are two-dimensional images in a one-dimensional parameter space, RCDT-POD can be applied to higher-dimensional data sets.

The results of our analysis show that RCDT-POD is a promising novel approach for compressing data sets of advection-dominated phenomena and travelling features and their subsequent predictive interpolation. 
The ability to capture and preserve, to some extent, certain geometrical features, within complex CFD data sets is, in fact, the most important property of this approach. For data sets with moving features (e.g., translation, scaling), this geometrical preservation is crucial.

The approach, however, introduces artefacts and errors due to the numerical implementation of the forward and inverse transform. We show that, while this intrinsic error in the discrete transform can be relatively large in the global metrics (e.g., $L^2$ norm), the POD reconstruction and interpolation error when predicting moving features can be orders of magnitude larger when linear ROM techniques are used. Furthermore, the artefacts and errors introduced by the discrete transforms can be reduced by improving the Radon and Cumulative Distribution discretisation and can be alleviated via  pre- and post-processing options. These will be considered in our future works.

% The results presented here on the physical interpolation error of a RCDT-POD workflow, when interpolating between variational flows of a chosen parameter, are enlightening on this aspect of the transform's use. Indicating a likelihood of interpolation error being the predominant error over the reconstructive and \textit{intrinsic} errors mentioned before.
%In future studies, the errors coming from the discrete transform, the POD truncation (reconstruction) and the interpolation will be further analysed separately to explore the RCDT-POD ROM applicability better when predicting advection-dominated flows of varying geometrical shapes, configurations, and speeds, and with more interpolated parameters.
Future work will focus on improving the RCDT accuracy for fields that do not go to zero at the boundary through filtering and thresholding. As suggested by the reviewers of this paper, which we gratefully acknowledge for their thorough and detailed work, this will include finding strategies alternative to the Radon transform to transform three-dimensional snapshots into one-dimensional distributions for the use of CDT. 
Another interesting aspect that deserves investigation is the extension using intrusive approaches. Such an extension would require an efficient way to compute the Jacobian of the non-linear RCDT transformation. In the online stage, the linear superposition of the spatial modes of the low-rank manifold will be characterised by a dependence on the non-linear map of the linear combination computed in the transformed space. When projecting the governing equation onto the reduced order manifold, the explicit form of the Jacobian will have to be computed, which is non-trivial (especially without sacrificing the computational efficiency).

%%%%%%%%%%%%%%%%%%%%%%%%%%%%%%
%%%%%%%%%%%%%%%%%%%%%%%%%%%%%%
\section*{Acknowledgements}
%%%%%%%%%%%%%%%%%%%%%%%%%%%%%%
%%%%%%%%%%%%%%%%%%%%%%%%%%%%%%
The authors thank Martina Cracco for her help in setting up the numerical algorithms. TL and RB have been funded
by the European Union through the project SILENTPROP, "Assessing noise generation in aircraft with distributed electric propulsion" (Grant agreement ID 882842). GS acknowledges the financial support under the National Recovery and Resilience Plan (NRRP), Mission 4, Component 2, Investment 1.1, Call for tender No. 1409 published on 14.9.2022 by the Italian Ministry of University and Research (MUR), funded by the European Union – NextGenerationEU– Project Title ROMEU – CUP P2022FEZS3 - Grant Assignment Decree No. 1379  adopted on 01/09/2023 by the Italian Ministry of Ministry of University and Research (MUR). The computations in this work have been performed with \textit{PyTransKit} \cite{abu_hasnat_mohammad_rubaiyat_pytranskit_2022} and \textit{EZyRB} \cite{demo_ezyrb_2018}; we acknowledge developers and contributors to both libraries.

\section*{Statements and Declarations}
The authors declared that they have no conflict of interest.

\section*{Author contributions}
TL: Writing - Original Draft, Software, Data processing. RB: Writing - Original Draft, Software, Data processing. RJL: Funding acquisition, Project administration. GS: Writing - Review \& Editing, Software, Methodology, Conceptualization, Supervision. MI: Writing - Review \& Editing, Software, Methodology, Conceptualization, Supervision, Project administration.

%%%%%%%%%%%%%%%%%%%%%%%%%%%%%%
%%%%%%%%%%%%%%%%%%%%%%%%%%%%%%
\bibliographystyle{unsrtnat}
\bibliography{bib/biblio, bib/paper_refs}

% % \newpage
% % %%%%%%%%%%%%%%%%%%%%%%%%%%%%%%
% % %%%%%%%%%%%%%%%%%%%%%%%%%%%%%%
% % %%%%%%%%%%%%%%%%%%%%%%%%%%%%%%
% % %%%%%%%%%%%%%%%%%%%%%%%%%%%%%%
% \begin{appendices}
%     \addtocontents{toc}{\protect\setcounter{tocdepth}{0}}
% % %%%%%%%%%%%%%%%%%%%%%%%%%%%%%%
% % %%%%%%%%%%%%%%%%%%%%%%%%%%%%%%
% % %%%%%%%%%%%%%%%%%%%%%%%%%%%%%%
% % %%%%%%%%%%%%%%%%%%%%%%%%%%%%%%

% % %%%%%%%%%%%%%%%%%%%%%%%%%%%%%%
% % %%%%%%%%%%%%%%%%%%%%%%%%%%%%%%
% \section{Additional Figures}
% %%%%%%%%%%%%%%%%%%%%%%%%%%%%%%
% %%%%%%%%%%%%%%%%%%%%%%%%%%%%%%

% \end{appendices}
% %%%%%%%%%%%%%%%%%%%%%%%%%%%%%%
% %%%%%%%%%%%%%%%%%%%%%%%%%%%%%%

\end{document}